\documentclass[12pt]{article}
\usepackage{latexsym}
\usepackage{amsmath}
\usepackage{amssymb}
\usepackage{amsthm}
\usepackage{mathtools}
\usepackage{dsfont}
\usepackage{bbm}
\usepackage{color}
\usepackage[tmargin=1.25in,bmargin=1.25in,lmargin=1.1in,rmargin=1.1in]{geometry}
\usepackage{tikz}
\usepackage{subcaption}
\usepackage{caption}
\usetikzlibrary{arrows,automata}
\usepackage{pstricks}
\usepackage{graphicx}
\usepackage{epsfig}
\usepackage{hyperref} 

\allowdisplaybreaks

\definecolor{Red}{rgb}{1,0,0}
\definecolor{halfgray}{gray}{0.55}
\definecolor{webgreen}{rgb}{0,.5,0}
\definecolor{webbrown}{rgb}{.6,0,0}
\definecolor{Maroon}{cmyk}{0, 0.87, 0.68, 0.32}
\definecolor{RoyalBlue}{cmyk}{1, 0.50, 0, 0}
\definecolor{Black}{cmyk}{0, 0, 0, 0}

\hypersetup{
    colorlinks=true, linktocpage=true, pdfstartpage=1, pdfstartview=FitV,
    breaklinks=true, pdfpagemode=UseNone, pageanchor=true, pdfpagemode=UseOutlines,
    plainpages=false, bookmarksnumbered, bookmarksopen=true, bookmarksopenlevel=1,
    hypertexnames=true, pdfhighlight=/O,
    urlcolor=webbrown, linkcolor=RoyalBlue, citecolor=webgreen,
}

\newtheorem{thm}{Theorem}[section]
\newtheorem{proposition}[thm]{Proposition}

\newtheorem{lem}[thm]{Lemma}
\newtheorem{definition}[thm]{Definition}
\newtheorem{remark}[thm]{Remark}

\makeatletter\@addtoreset{equation}{section}\makeatother

\begin{document}

\newcommand{\D}{\mathrm{d}}
\newcommand{\R}{\mathbb{R}}
\newcommand{\E}{\mathrm{E}}
\newcommand{\var}{\operatorname{Var}}
\newcommand{\cov}{\operatorname{Cov}}
\newcommand{\calG}{\mathcal{G}}
\newcommand{\calF}{\mathcal{F}}
\newcommand{\B}{\mathcal{B}}
\newcommand{\calC}{\mathcal{C}}
\newcommand{\calR}{\mathcal{R}}
\newcommand{\e}{\mathrm{e}}
\newcommand{\p}{\mathrm{P}}
\newcommand{\bbN}{\mathbb{N}}
\newcommand{\K}{\mathbb{K}}
\newcommand{\bbP}{\mathbb{P}}
\newcommand{\bbE}{\mathbb{E}}
\newcommand{\bbT}{\mathbb{T}}
\newcommand{\bbQ}{\mathbb{Q}}
\newcommand{\bfP}{\mathbf{P}}
\newcommand{\bfE}{\mathbf{E}}
\newcommand{\calE}{\mathcal E}
\newcommand{\GW}{\operatorname{GW}}
\newcommand{\Tra}{\operatorname{Tra}}
\newcommand{\RW}{\operatorname{RW}}
\newcommand{\RWP}{\RW\times\bfP}
\newcommand{\F}{\mathcal F}
\newcommand{\Bb}{\operatorname{Bb}}
\newcommand{\id}{\operatorname{id}}

\pagenumbering{arabic}

\newcommand{\Peaug}{\mathrm{P}_\varepsilon^{\mathrm{aug}}}
\newcommand{\Pehat}{\hat{\mathrm{P}}_\varepsilon}
\newcommand{\PPehat}{\hat{\mathbb{P}}_\varepsilon}
\newcommand{\Eehat}{\hat{\mathrm{E}}_\varepsilon}
\newcommand{\EEehat}{\hat{\mathbb{E}}_\varepsilon}


\title{Central limit theorem for a random walk on Galton-Watson trees with random conductances}
\author{Tabea Glatzel, Jan Nagel} 
\maketitle

\begin{abstract} 
We show a central limit theorem for random walk on a Galton-Watson tree, when the edges of the tree are assigned randomly uniformly elliptic conductances. When a positive fraction of edges is assigned a small conductance $\varepsilon$, we study the behavior of the limiting variance as $\varepsilon\to 0$. Provided that the tree formed by larger conductances is supercritical, the variance is nonvanishing as $\varepsilon\to 0$, which implies that the slowdown induced by the $\varepsilon$-edges is not too strong. The proof utilizes a specific regeneration structure, which leads to escape estimates uniform in $\varepsilon$.
\end{abstract}

{\bf Keywords:}
Functional central limit theorem, random walk in random environment, Galton-Watson trees

\smallskip

{\bf MSC 2020:} {
60K37, 
60F17, 
60K40 
}

\section{Introduction}
\label{sec:intro}

Random walks on Galton-Watson trees are a prominent example of random walks in random environments and are often viewed as a mean-field model for random walk on high-dimensional percolation clusters. When $T$ is a supercritical Galton-Watson tree conditioned to be infinite, the simple random walk on $T$ is transient and the distance to the root satisfies a law of large numbers with explicit and positive limiting velocity \cite{LyoPemPer96a}. A functional central limit theorem for the distance was proven by Piau \cite{Piau98CLT}. When the edges of the tree are randomly assigned i.i.d. edge weights, or conductances, the random walk in the environment formed by the tree and the conductances is a Markov chain which crosses an edge with probability proportional to the conductance of that edge. In this case, a law of large numbers still holds and was proven by \cite{gantert2012random}, where the limit cannot be computed explicitly even when the conductances can take only two different values. 

Such a random walk is particularly suited to model random processes in inhomogeneous environments. The potentially anomalous behavior of such random walks in random environments, for example anomalous scaling, aging and slowdown, has been well studied, for processes on trees see for example \cite{Dembo2002largeDeviations,Aidekon2010LargeDeviations,glatzel2021speed}. In particular we expect a trapping behaviour, when the conductances of some edges are very large or close to zero. The aim of this paper is to study the effect of small edge weights on the fluctuations of the random walk. 

We assign to the edges of a Galton-Watson tree randomly conductances, which are bounded and bounded away from 0. For the random walk in this environment, we show that the distance to the root satisfies a functional central limit theorem. The diffusivity $\sigma^2$ of the limit depends in an intricate way on the offspring law and the law of the conductances. In order to study this dependence, we assign a small conductance $\varepsilon>0$ to a positive fraction of edges. When $\varepsilon\to 0$, we show that $\sigma^2=\sigma^2(\varepsilon)$ is bounded away from zero, provided that the tree formed by edges with conductance larger than $\varepsilon$ is supercritical. This implies a qualitative estimate on the slowdown effect created by finite subtrees formed by 1-edges. In order to leave such a finite subtree, the random walk has to cross an $\varepsilon$-edge, which happens very rarely when $\varepsilon$ is small. Note that such a trapping phenomenon cannot occur when we reduce conductances to zero, which suggests a discontinuity of $\sigma^2(\varepsilon)$ at $\varepsilon=0$. In order to control the variance for small $\varepsilon$, we need a regeneration structure with moment bounds uniform in $\varepsilon$. The derivation of uniform bounds for escape probabilities and moments of regeneration times is the major challenge in the proofs. 

A closely related variant is the biased random walk on Galton-Watson trees. In this model, introduced by \cite{LyoPemPer96b}, edges between generation $n$ and $n+1$ are assigned conductance $\lambda^{-n}$, for a bias parameter $\lambda>0$. For small $\lambda$, this creates a drift away from the root, whereas for large $\lambda$, transition probabilities towards the root are larger and the random walk might become recurrent. A law of large numbers \cite{aidekon2014biased} and central limit theorem \cite{PerZei2008} hold also in this model, as well as trapping phenomenon. We refer to  \cite{arous2012biased,croydon2013biased,bowditch2018escape,bowditch2020differentiability}  or \cite{arous2016biased} for an overview. 

This paper is structured as follows. In the next section, we introduce the model more precisely and state our main results: the functional CLT in Theorem \ref{thm:CLT} and the estimate on the variance in Theorem \ref{thm:volatility}. In Section 3, we define the regeneration structure. The moment bounds for the regeneration times and distances are proven in Section 4, as well as some auxiliary estimates. The proof of the main results can be found in Section 5 and Section 6.

\section{Model and Results}
\label{sec:modelresults}

We denote by $\Omega$ be the set of all triples $(T,\rho,\xi)$, where $T$ is a tree with root $\rho$, $\mathcal{E}(T)$ its set of undirected edges and $\xi\in [0,\infty)^{\mathcal{E}(T)}$ is a configuration of conductances on the edges of $T$. An element $\omega\in \Omega$ is called an environment. 
We then let $\mathrm{P}$ be the law on $\Omega$, such that under $\mathrm{P}$, $T$ is a Galton-Watson tree with offspring law $\nu$ and root $\rho$, with an additional edge $(\rho,\rho^*)$ added to the root and conditioned on $T$, the conductances $\xi$ are independent and identically distributed with marginal law 
\begin{align} \label{def:conductancelaw}
\mu=\alpha\delta_\varepsilon+(1-\alpha)\mu_1
\end{align} 
for some $\alpha\in[0,1)$, where $\mu_1$ is a uniformly elliptic law with
\begin{align*}
\mu_1([\kappa^{-1},\kappa])=1
\end{align*}
for some ellipticity constant $\kappa\geq 1$ and such that $\mu_1(\{1\})>0$. We comment on this assumption after stating the main results. Without loss of generality, assume $\varepsilon<\kappa^{-1}$, so that we may refer to edges with ``small'' conductance $\varepsilon$ and ``large'' conductance in $[\kappa^{-1},\kappa]$.  

Let $N$ be the number of offspring of the root $\rho$. We assume that the tree is supercritical, $m_1=\mathrm{E}[N]>1$, it has no leaves, $\mathrm{P}(N=0)=0$ and all moments $m_k=\E[N^k]$ ($k>0$) exist. 
Note that the subtree formed by the edges with large conductance containing the root is supercritical if and only if $(1-\alpha)\E[N]>1$.
For vertices $u,v\in T$, we write $u\sim v$, if $u$ and $v$ are connected by an edge in $T$ and we say that $u$ and $v$ are neighbors. We call the vertices that lie on the path from the root to a vertex $v\in T$ the ancestors of $v$ and we denote by $v^*$ the ancestor that is also neighbor of $v$. Vertex $\rho^*$ is the ancestor of the root in the enlarged tree. Its addition to the Galton-Watson tree simplifies later notation, but is not relevant for the validity of limit theorems. Neighbors of vertex $v$, which are not ancestors are called descendants. We denote the graph distance of a vertex $v\in T,v\neq \rho^*$ to the root $\rho$ by $|v|$ and set $|\rho^*|=-1$. Then, $|v|$ is the generation of vertex $v$, with the convention that the generation of $\rho^*$ is $-1$. We denote the $k$-th generation of the tree $T$ by 
\begin{align*}
G_k\coloneqq G_k(T)\coloneqq\{v\in T(\rho):|v|=k\}
\end{align*}
and $G_{-1}\coloneqq \{\rho^*\}$ is the generation of the ancestor of the root $\rho$.

Given an environment $\omega=(T,\rho,\xi)$ and a vertex $v\in T$, we let $P^v_\omega$ be the quenched law of the Markov chain $(X_n)_{n\geq 0}$ on $T$, starting at $v$ and with transition probabilities 
\begin{align} \label{quenchedlaw}
P^v_\omega (X_{n+1}=z | X_n=u) = \frac{\xi(u,z)}{\sum_{w\sim u} \xi(u,w)}
\end{align}
if $u\sim z$. We write $P_\omega$ for the law of $(X_n)_{n\ge 0}$ starting at the root $\rho$. The annealed law $\mathbb{P}$ is then defined by 
\begin{align} \label{annealedlaw}
\mathbb{P}(A\times B) = \int_A P_\omega(B) \mathrm{d} \mathrm{P}(\omega)
\end{align}
for measurable subsets $A\times B$ of the set of all pairs $(\omega,T^{\mathbb{N}_0})$ (as usual, the $\sigma$-algebras are generated by finite-dimensional projections). 
As proven by \cite{gantert2012random}, we then have 
\begin{align} \label{defspeed}
\lim_{n\to \infty} \frac{|X_n|}{n} = v =v(\nu,\mu) \quad \mathbb{P}-\text{almost surely,}
\end{align}
where the speed $v(\nu,\mu)$ of the random walk depends in an intricate way on the offspring law and the conductance law. 

Our main result is a functional central limit theorem for the process
\begin{align}\label{def:processW}
W_t^n=\frac1{\sqrt n}\big(|X_{\lfloor tn\rfloor}| - \lfloor tn \rfloor v\big),\quad \text{for }t\in [0,1].
\end{align}
We consider  $(W_t^n)_{t\in [0,1]}$ as a random variable which takes values in the space $\mathbb D[0,1]$, endowed with the Skorokhod topology and its Borel $\sigma$-algebra.

\begin{thm}\label{thm:CLT}
There exists a constant $\sigma^2(\varepsilon)=\sigma^2(\varepsilon,\alpha,\mu_1)> 0$, such that 
\begin{align*}
\big(W_t^n\big)_{t\in [0,1]} 
\xrightarrow[n\to\infty]{d} \big( \sigma(\varepsilon) B_t\big)_{t\in [0,1]}
\end{align*}
under $\bbP$, with $\big( B_t\big)_{t\in [0,1]}$ a standard Brownian Motion.
\end{thm}

The following theorem investigates the behavior of $\sigma^2(\varepsilon)$ as $\varepsilon$ tends to zero. The variance is bounded away from 0 provided that the subtree $T_1(\rho)$ consisting of all vertices connected to $\rho$ by edges with conductance larger than $\varepsilon$ (determined by $\mu_1$) is supercritical, which is equivalent to $(1-\alpha)\E[N]>1$.

\begin{thm}\label{thm:volatility}
	If $\bbP(|T_1(\rho)|=\infty)>0$, we have for $\sigma^2(\varepsilon)$ from Theorem~\ref{thm:CLT}
	\begin{align*}
		\liminf_{\varepsilon\to 0} \sigma^2(\varepsilon) >0.
	\end{align*}
\end{thm}

\begin{remark}
\textup{It is natural to ask for the continuity of $\sigma^2$ in $\varepsilon=0$. When $\varepsilon=0$ (and $\alpha>0$), the conductances are not elliptic anymore and certain edges of the tree cannot be crossed by the random walk, in particular it could be restricted to a finite subtree. For this reason, we condition on $|T_1(\rho)|=\infty$ in the natural definition of the environment. This is the environment law in \cite{Piau98CLT} or \cite{gantert2012random}, but also in different models of random media, such as random walk on percolation clusters in $\mathbb{Z}^d$ \cite{berger2007quenched,mathieu2007quenched}. We may then define $\sigma^2(0)$ as the limit variance of the random walk in such a conditioned environment. Note that a CLT with variance $\sigma^2(0)$ has only been proven in the case $\mu_1=\delta_1$.  
Then, we do not expect that $\sigma^2(\varepsilon)\to \sigma^2(0)$ as $\varepsilon\to 0$. The reason is that for $\varepsilon$ small, the random walk will get trapped in finite subtrees with edge weights in $[\kappa^{-1},\kappa]$, surrounded by $\varepsilon$-edges. The probability of entering such a subtree is of order $\varepsilon$, but then the random walk needs time of order $1/\varepsilon$ to escape the trap. Thus, a slowdown created by such traps will be visible for arbitrarily small $\varepsilon$. For $\varepsilon=0$ however, these finite subtrees cannot be entered, which allows larger fluctuations of the process. We therefore expect that 
\begin{align*}
\liminf_{\varepsilon\to 0} \sigma^2(\varepsilon) < \sigma^2(0) . 
\end{align*}
A comparable statement has been proven for the speed $v(\varepsilon)$ of random walk on Galton-Watson trees with conductances distribution $\mu$ as in \eqref{def:conductancelaw} \cite{glatzel2021speed}. Therein, the convergence can be made more precise and we can show that
\begin{align*}
\lim_{\varepsilon\to 0} v(\varepsilon) = \beta v(0) . 
\end{align*} 
For a constant $\beta\in [0,1)$.} 
\end{remark}

\begin{remark}
\textup{Let us comment on the assumption $\mu(\{1\})>0$, which rules out a continuous conductance law. We use a positive density of $1$-edges to construct regeneration times, where after reaching a new generation via such an edge the random walk never backtracks. With this construction, the increments after the regeneration time are independent of the increments before, since the only edge influencing both has a fixed conductance. Several alternative definitions of regeneration times are possible, possibly leading to a weaker assumption. For example, as in \cite{shen2002asymptotic}, one could accept the dependence created by a single edge and obtain a Markov structure for the increments between regeneration times. Independent regeneration times are constructed by \cite{guo2016einstein} using a coin trick. Both approaches would be possible in our situation, but they create substantial technical difficulties when we want to quantify how moments of regeneration times depend on $\varepsilon$ close to 0. In order to keep the arguments more concise, we chose to add the assumption $\mu(\{1\})>0$. Note that while in general we work with a conductance law $\mu$ with two atoms in $\varepsilon$ and 1, if one is only interested in a functional CLT as in Theorem \ref{thm:CLT}, it is possible to set $\alpha=0$ such that $\mu$ has only one atom in 1. Finally, the choice for an atom in 1 is arbitrary, as any uniformly elliptic conductance distribution $\mu$ with at least one atom can be rescaled to satisfy $\mu(\{1\})>0$.} 
\end{remark}

\section{Regeneration structure}
\label{sec:regeneration}
We use the existence of a renewal structure to decouple the increments of a random walk in a random environment. This classical regeneration argument goes back to \cite{kesten1977renewal}. See also \cite{LyoPemPer96a} and \cite{Piau98CLT} for similar formulations for random walks on Galton-Watson trees. The advantage of the construction of regeneration times is that we obtain an i.i.d. sequence of increments between regeneration times, which allows us to apply classical limit theorems.
Since we consider random conductances, we need an additional condition on the local environment of the regeneration point to obtain independence. The regeneration times are constructed in the following way: We wait until the random walker reaches a generation by traversing an edge with conductance one for the first time. We call this the first potential regeneration time $\sigma_1$. If the walker never returns to the ancestor $X_{\sigma_1}^*$, we call $\sigma_1$ the first regeneration time $\tau_1$. Otherwise, if the walker returns to the previous generation at some time $R_1>\sigma_1$, $\sigma_1$ is not a regeneration time. Instead, we wait until the walker first reaches a new generation by crossing an edge with conductance one. We denote this second potential regeneration time by $\sigma_2$. If the walker never visits the ancestor $X_{\sigma_2}^*$ again, then $\sigma_2$ is our first regeneration time $\tau_1$. Otherwise, we repeat the above procedure. The transience of the random walk on the tree implies that the first regeneration time $\tau_1$ is almost surely finite. A repetition of this method yields an infinite sequence $(\tau_k)_{k\ge 1}$ of regeneration times. This construction guarantees that the walker sees disjoint parts of the tree between regeneration times, which is crucial to obtain independent increments. We will now define the regeneration times formally.

\begin{definition} \label{def:potentialregtimes}
Given an environment $\omega=(T,\rho,\xi)$ and a random walk $(X_n)_{n\ge 0}$ on $T$, we define a sequence of stopping times 
\begin{align*}
\sigma_1&=\inf\{n\ge 1: |X_n|>|X_m| \text{ for all } m<n, \, \xi(X_{n-1},X_n)=1\},\\
R_1&=\inf\{n\ge\sigma_1:|X_n|=|X_{\sigma_1}|-1\}
\end{align*}
and recursively for $k>1$
\begin{align*}
\sigma_k&=\inf\{n\ge R_{k-1}: |X_n|>|X_m| \text{ for all } m<n, \, \xi(X_{n-1},X_n)=1\},\\
R_k&=\inf\{n\ge\sigma_k:|X_n|=|X_{\sigma_k}|-1\}.
\end{align*}
We then have $\sigma_1\le R_1\le \sigma_2\le R_2\le\dots$
\end{definition}
We notice that this recursive construction yields finite stopping times until $R_k=\infty$ for the first time, meaning that the random walk never returns to the ancestor of the $k$-th potential regeneration point. As described above, this is the first regeneration of the random walk. Accordingly, we set
\begin{align}\label{def:regtime1}
K=\inf\{k\ge 1: R_k=\infty\},\quad \tau_1=\sigma_K.
\end{align}
It follows from Lemma \ref{lem:annealedescapeprob} that $K<\infty$ holds almost surely and $\tau_1$ is well-defined. We call $\tau_1$ the first regeneration time.
To define the $k$-th regeneration time, we let $\theta_m(v_n)_{n\ge 0} = (v_{n+m})_{n\ge 0}$ be the time shift of an infinite path $(v_n)_{n\ge0}$. For $k\ge 2$, the $k$-th regeneration time is then defined as
\begin{align*}
\tau_k=\tau_{k-1}+\tau_1\circ\theta_{\tau_{k-1}},
\end{align*}
where $\tau_1\circ\theta_{\tau_{k-1}}$ is the first regeneration time of the shifted path. Lemma~\ref{lem:existence} below shows that these regeneration times are well-defined. In particular, we obtain an infinite sequence $(\tau_k)_{k\ge 1}$ of regeneration times.

In order to define the regeneration structure as above, a uniform bound for the probability that the random walk hits the ancestor of the root is crucial. We denote this hitting time by 
\begin{align*}
\eta_*=\eta_{\rho^*}=\inf\{n\ge 0 : X_n=\rho^*\}.
\end{align*}
Moreover, given an environment $\omega=(T,\rho,\xi)$, we denote the subtree formed by a vertex $v\in T$ and all its descendants by $T(v)$ and the associated environment by $\omega(v)=\big(T(v),v,(\xi_e)_{e\in\mathcal E(T(v))}\big)$.
We write $T_1(v)$ for the subtree of $T(v)$ formed by the edges with conductance in $[\kappa^{-1},\kappa]$, that contains the vertex $v$. 

\begin{lem}[Annealed escape probability]\label{lem:annealedescapeprob}
There exists some constant $c_\varepsilon>0$, depending on $\nu$, $\alpha,\kappa$ and $\varepsilon$, such that
		\begin{align*}
			\bbP(\eta_*=\infty|\xi(\rho,\rho^*)=1)\ge c_\varepsilon .
		\end{align*}
If the subtree $T_1(\rho)$ is supercritical, i.e. if $(1-\alpha)\E[N]>1$, there exists some constant $c>0$, depending on $\nu, \kappa$ and $\alpha$, but independent of $\varepsilon$, such that
		\begin{align*}
			\bbP(\eta_*=\infty|\xi(\rho,\rho^*)=1)\ge c .
		\end{align*}
\end{lem}

\begin{figure}
	\centering
	\begin{tikzpicture}
	\draw[style=dashed] (0,0) to [out=60, in=180] (3.5,2.5)to [out=0, in=90] (6,0) to [out=270, in=0] (3.5,-2.5) to [out=180, in=300] (0,0);
	\draw (5.7,-2.1) node[] {$T(v)$};
	\draw[-, line width=0.5mm] (0,0) to (-2,-0.75);
	\draw[-, line width=0.5mm] (0,0) to (2,1);
	\draw[-, line width=0.5mm] (0,0) to (2,-1);
	\draw[-, line width=0.5mm] (2,1) to (4,1); \draw (4.5,1) node[] {$\dots$};
	\draw[-, line width=0.5mm] (2,1) to (4,1.7); \draw (4.5,1.7) node[] {$\dots$};
	\draw[-, line width=0.5mm] (2,1) to (4,0.3); \draw (4.5,0.3) node[] {$\dots$};
	\draw[-, line width=0.5mm] (2,-1) to (4,-1.5); \draw (4.5,-1.5) node[] {$\dots$};
	\draw[-, line width=0.5mm] (2,-1) to (4,-0.5); \draw (4.5,-0.5) node[] {$\dots$};
	\draw[-] (-2,-0.75) to (0,-1.5); 
	\draw[-] (-2,-0.75) to (-4,-2); 
	\draw (-0.1,-0.3) node[] {$v$};
	\draw (-2.1,-0.5) node[] {$v^*$};
	\draw (-1,0) node[] {$\boldsymbol{T^*(v)}$};
	\end{tikzpicture}
	\caption{$T(v)$ is the subtree formed by $v$ and all its descendants, $T^*(v)$ (indicated by the thick edges) is the subtree composed of $T(v)$ and $v^*$.}
	\label{fig:subtrees}
\end{figure}
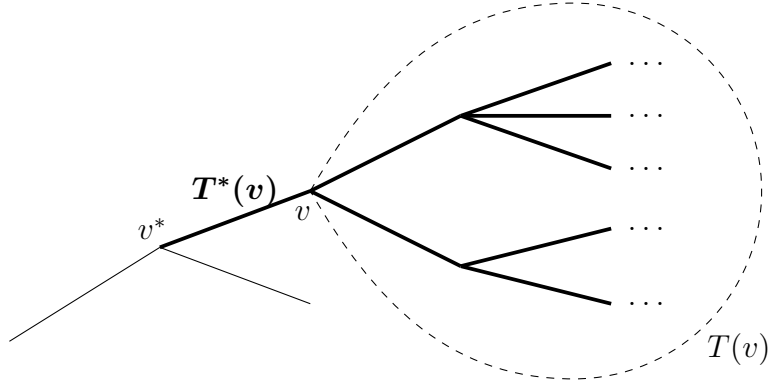

To formulate the decoupling property of regeneration times, we have to introduce some more notations. Given an environment $\omega=(T,\rho,\xi)$, we let $T^*(v)$ be the subtree composed of $T(v)$ and the ancestor $v^*$ of $v$. Since $\rho^*$ does not have an ancestor, we set $T^*(\rho^*)=T(\rho^*)$. The associated environment is denoted by $\omega^*(v) = \big(T^*(v),v,(\xi_e)_{e\in\mathcal E(T^*(v))}\big)$ (see Figure~\ref{fig:subtrees}).
To encode the position of a vertex in this subtree, we think of a tree as a subset of an extended Ulam-Harris tree $\bbT^*=\bbT\cup\{-1\}$, where $-1$ denotes the artificial ancestor of the root. This means that each vertex corresponds to a sequence of integers, we may identify the subtree $T^*(v)$, rooted at $v$, with the set 
\begin{align*}
	[T^*(v)] = [(T^*(v),v)] = \{z-v : z\in T^*(v)\} \subset \bbT^*.
\end{align*}
Here, the difference between two vertices $z=(z_1,\dots,z_n)$ and $v=(z_1,\dots,z_l)$ is defined by
\begin{align*}
	z-v = \begin{cases}
		(z_{l+1},\dots, z_n), & z\in T(v)\\
		-1, & z=v^*.
	\end{cases}
\end{align*}
In other words, $[T^*(v)]$ shifts the subtree $T^*(v)$ such that the root $v$ is identified with $\emptyset$ in the extended Ulam-Harris tree. Therefore, the structure of the tree remains unchanged, but the tree becomes independent of the position of the vertex $v$ in the original tree. We write $[\omega^*(v)]$ for corresponding environment.

Further, we denote by $\calG$ the $\sigma$-field on the set of environments $\Omega$ and by $\calF$ the $\sigma$-field on the set of trajectories $\mathbb{T}^{\mathbb{N}_0}$. For a set $B\in\calG\otimes\calF$ we define the time shift $\Theta_k$ by
\begin{align*}
B\circ\Theta_{\tau_k}=\{(\omega,(X_n)_{n\ge 0}):([\omega^*(X_{\tau_k})], (X_n-X_{\tau_k})_{n\ge \tau_k})\in B\}.
\end{align*}
We note that $X_n-X_{\tau_k}$, which denotes the position of $X_n$ in the subtree $[T^*(X_{\tau_k})]$, is well-defined by the definition of $\tau_k$. 
Moreover, for $k\in\bbN$ we introduce the $\sigma$-field $\calG_k$, which is generated by the sets
\begin{align*}
\{\tau_k=m, X_0=v_0,\dots,X_m=v_m, \omega\backslash \omega(v_m)\in A\}.
\end{align*}
Here, for two environments $\omega_1=(T_1,\rho_1,(\xi_e)_{e\in\calE(T_1)}),\,\omega_2=(T_2,\rho_2,(\xi_e)_{e\in\calE(T_2)})$, such that $T_2$ is a subtree of $T_1$, we write $T_1\backslash T_2$ for the tree, which we obtain when we remove from $T_1$ all edges of $T_2$ and all isolated vertices. The associated environment is denoted by $\omega_1\backslash\omega_2 = (T_1\backslash T_2, \rho_1, (\xi_e)_{e\in\calE(T_1\backslash T_2)})$. 

The following lemmas give the existence of an infinite sequence of regeneration times $(\tau_k)_{k\ge 1}$ and the independence of the inter-regeneration epochs. Their proofs are standard and follow as in \cite[Lemma 17 and Proposition 18]{guo2016einstein}.

\begin{lem}[Existence of regeneration times]\label{lem:existence}
For any $k\ge 1$, the $k$-th regeneration time $\tau_k$ is almost surely finite. Moreover, we have
\begin{align*}
\bbP(B\circ \Theta_{\tau_k}|\calG_k)=\bbP(B|\eta_*=\infty,\xi(\rho,\rho^*)=1)
\end{align*}
for $B\in\calG\otimes\calF$.
\end{lem}


\begin{proposition}[Stationarity and independence]\label{prop:stationaryindpendent}
Under $\bbP$, the sequence
\begin{align*}
\big([\omega^*(X_{\tau_n})\backslash\omega(X_{\tau_{n+1}})], (X_k-X_{\tau_n})_{\tau_n\le k\le \tau_{n+1}}, \tau_{n+1}-\tau_n\big)_{n\ge 1}
\end{align*}
is stationary and independent. 
Furthermore, the marginal distribution of this sequence is given by
\begin{align*}
&\bbP\big( [\omega^*(X_{\tau_n})\backslash\omega(X_{\tau_{n+1}})]\in A_1, (X_k-X_{\tau_n})_{\tau_n\le k\le \tau_{n+1}}\in A_2, \tau_{n+1}-\tau_n\in A_3\big)\\
&\quad= \bbP\big(\omega\backslash\omega(X_{\tau_1})\in A_1, (X_k)_{k\le \tau_1}\in A_2, \tau_1\in A_3\big|\eta_*=\infty,\xi(\rho,\rho^*)=1\big).
\end{align*}
\end{proposition}

The regeneration structure implies a central limit theorem, provided that the regeneration times and the regeneration distances are sufficiently integrable. The necessary moment bounds are given in the next two lemmas.

\begin{lem}[Moment bounds on regeneration distances]\label{lem:momentboundregdist}
For $\varepsilon>0$ fixed and for any $q\ge 1$, there exists a constant $C(\varepsilon)>0$ depending on $\varepsilon$, such that
	\begin{align*}
		\bbE\big[|X_{\tau_1}|^q\big]\le C(\varepsilon) \quad \text{and} \quad \bbE\big[(|X_{\tau_2}|-|X_{\tau_1}|)^q\big]\le C(\varepsilon).
	\end{align*}
	If the subtree $T_1(\rho)$ is supercritical, i.e. if $\p(|T_1(\rho)|=\infty)>0$, the bounds above hold for a constant $C$ independent of $\varepsilon$. 
\end{lem}

\begin{lem}[Moment bounds on regeneration times]\label{lem:momentboundregtimes}
	For $\varepsilon>0$ fixed and for any $q\ge1$ there exists a constant $C(\varepsilon)>0$, depending on $\varepsilon$, such that
	\begin{align*}
	\bbE\big[\tau_1^q\big]\le C(\varepsilon) \quad \text{and} \quad \bbE\big[(\tau_2-\tau_1)^q\big]\le C(\varepsilon).
	\end{align*}
\end{lem}

For the proof of Lemma~\ref{lem:momentboundregdist} we need a quenched bound for the probability that the random walk, starting at some vertex in generation $L$, hits the ancestor of the root at some time. Unfortunately, a quenched version of Lemma~\ref{lem:annealedescapeprob} does not hold. Instead, we can bound the expected number of vertices $v$ in the $L$-th generation that have a small probability to never hit the ancestor of the root. 

\begin{lem}[quenched indirect escape probability] \label{lem:quenchedindirectescapeprob}  
There exists a function $\gamma$ with $\gamma(L,\nu,\alpha)\to 0 $ as $L\to\infty$ and a constant $C=C(\nu,\alpha,\kappa)>0$, both independent of $\varepsilon$, such that for $0<\delta<\frac1{1+\kappa^2L}$
\begin{align*}
\E\Bigg[ \sum_{v\in G_L}\mathds 1_{\{P^{v}_{\omega\backslash \omega(v)}(\eta_*=\infty)< \delta\}} \Bigg]
&\le 
\gamma(L,\nu,\alpha) + C\frac{\delta\kappa}{1-\delta(1+\kappa^2L)} m^{L}.
\end{align*}
\end{lem}

To proof Lemma~\ref{lem:momentboundregtimes}, we need a bound on the probability that the random walk returns $n$ times to the root. We denote by  
\begin{align} \label{def:localtime}
 	L(z) = \sum_{k=0}^\infty \mathds 1_{\{X_k=z\}}
\end{align} 
the local time of a vertex $z$. The required estimate is given in the next lemma.
\begin{lem}[bound on local time]\label{lem:returnprob}
	For any $\beta >0$ there exists some constant $C_\varepsilon>0$ such that
	\begin{align*}
		\bbP(L(\rho)>n|\xi(\rho,\rho^*)=\kappa)
		\le C_\varepsilon n^{-\beta}.
	\end{align*}
\end{lem}

\section{Proofs of auxiliary results}
\label{sec:technicalproofs}

\subsection{Proof of Lemma \ref{lem:annealedescapeprob}}

Let us assume that the tree $T_1(\rho)$ is supercritical, i.e. $(1-\alpha)\E[N]>1$ and the probability that $T_1(\rho)$ survives is strictly positive. Thus, we can bound the annealed escape probability as follows
\begin{align}\label{hittingprob1}
\bbP(\eta_*=\infty|\xi(\rho,\rho^*)=1) \notag
&\ge \bbP\big(\eta_*=\infty, |T_1(\rho)| =\infty \big|\xi(\rho,\rho^*)=1\big)\\ 
&=\bbP\big(\eta_*=\infty\big|\xi(\rho,\rho^*)=1, |T_1(\rho)|=\infty\big) \p(|T_1(\rho)|=\infty),
\end{align}
where we used the independence of the conductance of the edge $(\rho,\rho^*)$ and the subtree $T_1(\rho)$. 
Recalling the definition of the annealed law $\bbP$ in \eqref{annealedlaw}, we can write
\begin{align}\label{hittingprob2}
\bbP(\eta_*<\infty|\xi(\rho,\rho^*)=1, |T_1(\rho)|=\infty) 
=\E[P_\omega(\eta_*<\infty)|\xi(\rho,\rho^*)=1, |T_1(\rho)|=\infty].
\end{align}
The quenched return probability in \eqref{hittingprob2} is the limit of $P_\omega(\eta_*<\eta_k)$ as $k\to \infty$, 
where $\eta_k\coloneqq \inf\{n\ge 0:|X_n|=k\}$ is the hitting time of the $k$-th generation of the tree.

Let $\calC_\omega(v,A)$ be the effective conductance between a vertex $v\in T$ and a set of vertices $A\subset T$. The associated effective resistance is denoted by $\calR_\omega(v,A)=\calC_\omega(v,A)^{-1}$. For a precise definition and the relation between random walks and electric networks, we refer to \cite{doyle1984random} and \cite{LyoPer16}. If $\xi(\rho,\rho^*)=1$, we have
\begin{align*}
P_\omega(\eta_*<\eta_k) = \frac{\calC_\omega(\rho,\rho^*)}{\calC_\omega(\rho,\rho^*)+\calC_\omega(\rho,G_k)} = \frac1{1+\calC_\omega(\rho,G_k)}.
\end{align*}
Let $T_1^{\Bb}=T_1^{\Bb}(\rho)$ denote the backbone tree of $T_1(\rho)$, where all vertices that do not have an infinite line of descent are removed. We write $\omega_1^{\Bb}=\big(T_1^{\Bb},\rho,(\xi_e)_{e\in\calE(T_1^{\Bb})}\big)$ for the associated environment. Rayleigh's monotonicity principle implies
\begin{align*}
\calC_\omega(\rho,G_k)\ge \calC_{\omega_1^{\Bb}}(\rho,G_{k}).
\end{align*}
Hence, we have
\begin{align*}
P_\omega(\eta_*<\eta_k) 
\le \frac1{1+\calC_{\omega_1^{\Bb}}(\rho,G_k)}
= \frac{\calR_{\omega_1^{\Bb}}(\rho,G_{k})}{1+\calR_{\omega_1^{\Bb}}(\rho,G_{k})}
\end{align*}
and
\begin{align*}
P_\omega(\eta_*<\infty) 
\le  \frac{\calR_{\omega_1^{\Bb}}(\rho,\infty)}{1+\calR_{\omega_1^{\Bb}}(\rho,\infty)}.
\end{align*}
The effective resistance $\calR_{\omega_1^{\Bb}}(\rho,\infty)$ only depends on the subtree rooted at $\rho$. Thus, it is independent of the edge $(\rho,\rho^*)$ under $\p$. Together with \eqref{hittingprob2}, this implies
\begin{align*}
\bbP(\eta_*<\infty|\xi(\rho,\rho^*)=1, |T_1(\rho)|=\infty) 
&\le \E\left[\frac{\calR_{\omega_1^{\Bb}}(\rho,\infty)}{1+\calR_{\omega_1^{\Bb}}(\rho,\infty)}\middle|\xi(\rho,\rho^*)=1, |T_1(\rho)|=\infty  \right]\\
&= \E\left[\frac{\calR_{\omega_1^{\Bb}}(\rho,\infty)}{1+\calR_{\omega_1^{\Bb}}(\rho,\infty)}\middle||T_1(\rho)|=\infty  \right].
\end{align*}
Conditioned on the the survival of $T_1(\rho)$, the backbone tree $T_1^{\Bb}$ is a supercritical Galton-Watson tree without leaves and with resistances bounded by $\kappa$ (see Proposition~4.10 in \cite{Lyon92}). Hence, Lemma~9.1 in \cite{LyoPemPer96a} implies 
\begin{align} \label{resistancemoment1tree}
\E[\calR_{\omega_1^{\Bb}}(\rho,\infty)||T_1(\rho)|=\infty ]\le C
\end{align}
for some constant $C=C(\nu,\alpha,\kappa)<\infty$. We obtain with Jensen's inequality
\begin{align*}
\bbP(\eta_*<\infty|\xi(\rho,\rho^*)=1, |T_1(\rho)|=\infty) 
&\le \frac{\E[\calR_{\omega_1^{\Bb}}(\rho,\infty)||T_1(\rho)|=\infty]}{1+\E[\calR_{\omega_1^{\Bb}}(\rho,\infty)||T_1(\rho)|=\infty]}\le\frac C{1+C} < 1.
\end{align*}
In total, due to \eqref{hittingprob1}, we have
\begin{align*}
\bbP(\eta_*=\infty|\xi(\rho,\rho^*)=1)\ge \left(1-\frac C{1+C}\right) \p(|T_1(\rho)|=\infty) > 0.
\end{align*}

If the  tree $T_1(\rho)$ is critical or subcritical, i.e. $(1-\alpha)\E[N]\leq1$, then the same arguments hold with the constant $C$ in \eqref{resistancemoment1tree} replaced by $C(\varepsilon)$, the (finite) expected resistance of the tree $T(\rho)$.  
\hfill $\Box$

\subsection{Proof of Lemma~\ref{lem:momentboundregdist}}
We prove the statement under the assumption that the subtree $T_1(\rho)$ is supercritical, that is, $T_1(\rho)$ has positive probability to survive. 
If the subtree $T_1(\rho)$ dies out with probability 1, the proof works analogously but the moment bounds depend on $\varepsilon$, since the lower bound in Lemma~\ref{lem:annealedescapeprob} depends on $\varepsilon$ in this case.

To show the moment bound on the regeneration distances, we follow the arguments of \cite{SziZer99LLN} with modifications as in \cite{HHN2020random}. 
By Lemma~\ref{lem:annealedescapeprob} and Proposition~\ref{prop:stationaryindpendent}, we have
\begin{align*}
\bbE\big[(|X_{\tau_2}|-|X_{\tau_1}|)^q\big]
&=\bbE\big[(|X_{\tau_1}|-|X_0|)^q\big|\eta_*=\infty,\xi(\rho,\rho^*)=1\big]\\
&\le \bbE\big[(|X_{\tau_1}|-|X_0|)^q\big] \bbP(\eta_*=\infty|\xi(\rho,\rho^*)=1)^{-1} \bbP(\xi(\rho,\rho^*)=1)^{-1}\\
&\le C \bbE\big[|X_{\tau_1}|^q\big].
\end{align*}
Hence, it is sufficient to bound $\bbE\big[|X_{\tau_1}|^q\big]$. Since the first regeneration time is the last finite potential regeneration time (see definition~\eqref{def:regtime1}), we can write 
\begin{align*}
\bbE\big[|X_{\tau_1}|^q\big]
= \sum_{k\in\bbN} \bbE\big[|X_{\tau_1}|^q\mathds 1_{\{\sigma_k<\infty,R_k=\infty\}}\big]
\le \sum_{k\in\bbN} \bbE\big[|X_{\sigma_k}|^q\mathds 1_{\{\sigma_k<\infty\}}\big].
\end{align*}
Applying the Cauchy-Schwarz inequality gives rise to the following upper bound
\begin{align}\label{momentbountdist-1}
\bbE\big[|X_{\tau_1}|^q\big]
\le \sum_{k\in\bbN} \bbE\big[|X_{\sigma_k}|^{2q}\mathds 1_{\{\sigma_k<\infty\}}\big]^{\frac12}\bbP(\sigma_k<\infty)^{\frac12}.
\end{align}
In order that $\sigma_k$ is finite, the random walker has to see a new generation and return to its ancestor at least $k-1$ times. Hence, Lemma~\ref{lem:annealedescapeprob} implies
\begin{align}\label{prob_sigmakfinite}
\bbP(\sigma_k<\infty) 
\le \bbP(\eta_*<\infty|\xi(\rho,\rho^*)=1)^{k-1} 
\le (1-c)^{k-1}
\end{align}
for a constant $c\in (0,1)$. To bound the expectation in \eqref{momentbountdist-1} we decompose the trajectory of the random walk. We introduce
\begin{align*}
M_k&\coloneqq|X_{\sigma_k}|,\\
N_k&\coloneqq\max\{m:\eta_m<R_k\}.
\end{align*}
That is, $M_k$ denotes the distance of the $k$-th potential regeneration point to the root. $N_k$ is the maximal generation visited by the random walk until it returns to the ancestor of $X_{\sigma_k}$. Further, we set
\begin{align*}
H_1\coloneqq M_1,\quad 
H_k\coloneqq M_k-N_{k-1},
\end{align*}
so that $H_k$ indicates the number of generations that are visited by random walk after exceeding the previous maximal generation $N_{k-1}$ and until reaching the next potential regeneration point. Lastly, we define
\begin{align*}
\tilde N_0\coloneqq 0,\quad 
\tilde N_k\coloneqq N_k-M_k.
\end{align*}
That is, $\tilde N_k$ counts the number of generations that are reached by the random walk for the first time after it exceeds the $k$-th potential regeneration point and until it returns to its ancestor. In other words, $\tilde N_k$ indicates the number of new generations that are visited during this excursion. Now, we can express the distance of the $k$-th regeneration point to the root by 
\begin{align*}
|X_{\sigma_k}|=M_k = \sum_{i=0}^{k-1} (\tilde N_i+H_{i+1}).
\end{align*}
Using the upper bounds in \eqref{momentbountdist-1} and \eqref{prob_sigmakfinite}, we obtain
\begin{align*}
\bbE\big[|X_{\tau_1}|^q\big]
\le \sum_{k\in\bbN} \bbE\left[\left(\sum_{i=0}^{k-1} (\tilde N_i+H_{i+1})\right)^{2q}\mathds 1_{\{\sigma_k<\infty\}}\right]^{\frac12}(1-c)^{\frac{k-1}2}.
\end{align*}
Applying Jensen's inequality twice gives the following upper bound
\begin{align*}
\bbE\big[|X_{\tau_1}|^q\big]
&\le \sum_{k\in\bbN} \bbE\left[ k^{2q-1}\sum_{i=0}^{k-1} \big(\tilde N_i+H_{i+1}\big)^{2q} \mathds 1_{\{\sigma_k<\infty\}}\right]^{\frac12} (1-c)^{\frac{k-1}2}\\
&\le \sum_{k\in\bbN} \bbE\left[ (2k)^{2q-1}\sum_{i=0}^{k-1} \big(\tilde N_i^{2q}+H_{i+1}^{2q}\big) \mathds 1_{\{\sigma_k<\infty\}}\right]^{\frac12} (1-c)^{\frac{k-1}2}.
\end{align*}
Noting that $\{\sigma_k<\infty\}\subseteq\{R_{k-1}<\infty\}\subseteq\{R_i<\infty\}$ holds for any $i\le k-1$, we arrive at
\begin{align}\label{momentbountdist-2}
\bbE\big[|X_{\tau_1}|^q\big]
\le \sum_{k\in\bbN} (2k)^{\frac{2q-1}2}\left(\sum_{i=0}^{k-1} \bbE\big[\tilde N_i^{2q}\mathds 1_{\{R_i<\infty\}}\big] + \sum_{i=0}^{k-1} \bbE\big[H_{i+1}^{2q}\mathds 1_{\{R_i<\infty\}}\big] \right)^{\frac12} (1-c)^{\frac{k-1}2}.
\end{align}
The right hand side remains finite, if we can uniformly bound the occurring moments.

\subsubsection{Uniform bound for $\bbE\big[H_{i+1}^{2q}\mathds 1_{\{R_i<\infty\}}\big]$}
Recall that $H_i$ counts the number of new generations visited by the walker after surpassing the previous maximal generation, until reaching the next potential regeneration point. In other words, after reaching the previous maximal generation, we wait until the random walk gets to a new generation by crossing an $1$-edge for the first time. 

Denote by 
\begin{align*}
V_1\coloneqq V_1(T) \coloneqq \{v\in T: \xi(v,z)=1\text{ for all } z\in G_1(T(v)))\}
\end{align*}
the set of all vertices $v\in T$ where each edge that connects $v$ with a descendant has conductance one. If the random walk reaches a new generation at vertex $v$, there is a positive probability that $v\in V_1$. In this case, there is a positive lower bound independent of $\varepsilon$ for the probability that the random walker moves to a descendant of $v$. If both of these events occur, the hitting time of generation $|v|+1$ is a potential regeneration. 

Therefore, it is possible to stochastically dominate 
\begin{align*}
H_{i+1}\mathds 1_{\{R_i<\infty\}}
\end{align*}
by a geometric random variable with success probability $p>0$ independent of $\varepsilon$, which implies 
\begin{align}\label{momentbound-H}
\bbE[(H_{i+1}\mathds 1_{\{R_i<\infty\}})^{2q}] 
\le C_1<\infty
\end{align}
for some constant $C_1 \ge 0$.

\subsubsection{Uniform bound for $\bbE\big[\tilde N_i^{2q}\mathds 1_{\{R_i<\infty\}}\big]$}
In view of \eqref{momentbountdist-2}, the proof is completed once we bound $\bbE[\tilde N_i^{2q}\mathds 1_{\{R_i<\infty\}}]$ uniformly in $i$.
Recall that $\tilde N_i$ counts the number of generations that are visited during the excursions starting from the potential regeneration points. From the construction of the potential regeneration points follows that these excursions  occur in disjoint parts of the tree. Hence, the random variables $\tilde N_i \mathds 1_{\{R_i<\infty\}}$ are i.i.d. and it suffices to show
\begin{align*}
\bbE\big[\tilde N_1^{2q}\mathds 1_{\{R_1<\infty\}}\big]
\le  c_2<\infty.
\end{align*}
This is implied by a bound for 
\begin{align}\label{moments_exptildeN}
\bbE\big[\e^{s\tilde N_1}\mathds 1_{\{R_1<\infty\}}\big]
=\sum_{n\in\bbN} \e^{sn}\bbP(\tilde N_1=n,R_1<\infty)
\end{align}
for some $s>0$. We show that the event $\{\tilde N_1=n,R_1<\infty\}$ has exponentially small probability.

Set
\begin{align*}
T\coloneqq \max\{m\ge 1:\eta_m<\eta_{X_0^*}\}-|X_0|,
\end{align*}
so that $T$ counts the number of generations visited by the random walk until it visits the ancestor of its starting point for the first time. 
Then, we calculate
\begin{align*}
&\bbP(\tilde N_1=n,R_1<\infty)
= \sum_{v\in\bbT} \E\big[P_\omega(\tilde N_1=n, R_1<\infty, \sigma_1<\infty, X_{\sigma_1}=v) \mathds 1_{\{\xi(v,v^*)=1\}}\big]\\
&=\sum_{v\in\bbT} \E\big[P^v_{\omega^*(v)}(T=n, \eta_{v^*} <\infty) P_\omega(\sigma_1<\infty,X_{\sigma_1}=v)\mathds 1_{\{\xi(v,v^*)=1\}}\big]\\
&=\sum_{v\in\bbT} \E\big[P^v_{\omega^*(v)}(T=n, \eta_{v^*} <\infty) P_\omega(\sigma_1<\infty,X_{\sigma_1}=v)\big| \xi(v,v^*)=1\big] \p(\xi(v,v^*)=1),
\end{align*}
where the second equality holds because of the strong Markov property. Since the random variables $P^v_{\omega^*(v)}(T=n, \eta_{v^*} <\infty)$ and $P(\sigma_1<\infty,X_{\sigma_1}=v)$ are independent under $\p(\,\cdot\,|\xi(v,v^*)=1)$, we get
\begin{align*}
&\bbP(\tilde N_1=n,R_1<\infty)\\
&=\sum_{v\in\bbT} \E\big[P^v_{\omega^*(v)}(T=n, \eta_{v^*} <\infty)\big| \xi(v,v^*)=1\big] \E\big[P_\omega(\sigma_1<\infty,X_{\sigma_1}=v)\mathds 1_{\{\xi(v,v^*)=1\}}\big]\\
&=\E\big[P_\omega(T=n, \eta_* <\infty)\big| \xi(\rho,\rho^*)=1\big] \bbP(\sigma_1<\infty)\\
&=\E\Bigg[\sum_{v\in G_{n}(T)} P_\omega(\eta_{n}<\eta_*<\eta_{n+1}, X_{\eta_{n}}=v)\Bigg| \xi(\rho,\rho^*)=1\Bigg]\\
&\le \E\Bigg[\sum_{v\in G_{n}(T)} P_\omega(\eta_*\circ\theta_{\eta_{n}} <\infty, X_{\eta_{n}}=v)\Bigg| \xi(\rho,\rho^*)=1\Bigg],
\end{align*}
where $\eta_*\circ\theta_{\eta_{n}}=\inf\{m\ge\eta_{n}:X_m=\rho^*\}$ denotes the hitting time of $\rho^*$ of the shifted path.
Using the strong Markov property in $\eta_{n}$, we arrive at
\begin{align}\label{bound-backtrack}
\bbP(\tilde N_1=n,R_1<\infty)
& \le\E\Bigg[\sum_{v\in G_{n}(T)} P^v_\omega(\eta_*<\infty) P_\omega( X_{\eta_{n}}=v) \Bigg| \xi(\rho,\rho^*)=1\Bigg].
\end{align}

\subsubsection{Quenched return probability}
We need a quenched bound for the return probability $ P^v_\omega(\eta_*<\infty)$. Fix a vertex $v\in G_{n}(T)$ in the $n$-th generation, set $r=\lfloor\frac n L\rfloor$ and denote the ancestor of $v$ in $G_{[i]}=G_{iL}(T)$ by $v_i$ for $i=0,\dots,r$. Further, we define 
\begin{align}
\eta_z =\inf\{n\ge 0: X_n=z\}
\end{align}
to be the hitting time of the vertex $z\in T$. Then, 
\begin{align}\label{quenchedreturnprob-1}
P^v_\omega(\eta_*<\infty) 
\le \prod_{i=1}^{r} P^{v_{i}}_\omega(\eta_{v_{i-1}}<\infty).
\end{align}
Using  Rayleigh's monotonicity principle, we have
\begin{align*}
P^{v_{i}}_\omega(\eta_{v_{i-1}}<\infty)
\le P^{v_{i}}_{\omega\backslash\omega(v_{i})}(\eta_{v_{i-1}}<\infty),
\end{align*}
such that
\begin{align*}
P^{v}_\omega(\eta_*<\infty)
\le \prod_{i=1}^{r} P^{v_{i}}_{\omega\backslash\omega(v_{i})}(\eta_{v_{i-1}}<\infty).
\end{align*}
We can establish an exponentially small bound for the probability above, if the underlying environment is good enough in some sense. Roughly speaking, we call an environment good, if for all vertices $v\in G_{[r]}$ the fraction of ancestors $v_i$ with a sufficiently high escape probability is large enough. Recall that $v_i$ denotes the ancestor of $v$ in generation $G_{[i]}$. More precisely, we introduce the set of good environments
\begin{align*}
B_r(\delta,\beta, L)
=\left\{ \sum_{i=1}^r \mathds 1_{\big\{P^{v_i}_{\omega\backslash\omega(v_i)} (\eta_{v_{i-1}}=\infty) \ge \delta \big\}}\ge\beta r \text{ for all } v\in G_{[r]}\right\}.
\end{align*}
Then, for $\omega\in B_r(\delta,\beta, L)$, any $v\in G_{n}(T)$ has at least $\beta r$ ancestors $v_i$, where $r=\lfloor\frac n L \rfloor$, such that the indirect escape probability $P^{v_{i}}_{\omega\backslash\omega(v_{i})}(\eta_{v_{i-1}}=\infty)$ is at least $\delta$. This implies
\begin{align*}
P^{v}_\omega(\eta_*<\infty)
&\le \prod_{i=1}^{r} \left( (1-\delta)\mathds 1_{\big\{P^{v_i}_{\omega\backslash\omega(v_i)} (\eta_{v_{i-1}}=\infty) \ge \delta \big\}} + \mathds 1_{\big\{P^{v_i}_{\omega\backslash\omega(v_i)} (\eta_{v_{i-1}}=\infty) <\delta \big\}} \right)\\
&\le (1-\delta)^{\beta r}.
\end{align*}
Together with the upper bound in \eqref{bound-backtrack}, we arrive at
\begin{align}\label{bound-backtrack-2} \notag
&\bbP(\tilde N_1=n,R_1<\infty)\\ \notag
& \le (1-\delta)^{\beta r} \E\Bigg[\sum_{v\in G_{n}(T)} P_\omega( X_{\eta_{n}}=v) \Bigg| \{\xi(\rho,\rho^*)=1\}\cap B_r(\delta,\beta, L) \Bigg]\\ \notag
&\quad+ \E\Bigg[\sum_{v\in G_{n}(T)} P_\omega( X_{\eta_{n}}=v) \Bigg| \{\xi(\rho,\rho^*)=1\}\cap B_r(\delta,\beta, L)^c \Bigg] \p(B_r(\delta,\beta, L)^c|\xi(\rho,\rho^*)=1)\\\notag
& = (1-\delta)^{\beta r} \E\big[P_\omega( \eta_{n}<\infty) \big| \{\xi(\rho,\rho^*)=1\}\cap B_r(\delta,\beta, L) \big]\\ \notag
&\quad+ \E\big[P_\omega( \eta_{n}<\infty) \big| \{\xi(\rho,\rho^*)=1\}\cap B_r(\delta,\beta, L)^c \big] \p(B_r(\delta,\beta, L)^c|\xi(\rho,\rho^*)=1)\\ 
& \le (1-\delta)^{\beta r} 
+ \p(B_r(\delta,\beta, L)^c|\xi(\rho,\rho^*)=1).
\end{align}
Having in mind that our aim is to show that the probability $\bbP(\tilde N_1=n,R_1<\infty)$ decays exponentially (see \eqref{moments_exptildeN}), we are done once we show that $B_r(\delta,\beta, L)^c$ has exponentially small probability.

\subsubsection{$B_r(\delta,\beta, L)^c$ has exponentially small probability}
To bound the probability of $B_r(\delta,\beta, L)^c$, we use similar arguments as \cite{grimmett1984random} and \cite{Dembo2002largeDeviations}. 
We introduce for $v\in G_{[r]}$
\begin{align*}
A_r(\delta,v, L) \coloneqq \sum_{i=1}^r \mathds 1_{\{P^{v_i}_{\omega\backslash \omega(v_i)}(\eta_{v_{i-1}}=\infty)\ge \delta)\}}
\end{align*}
and
\begin{align*}
Z_r(\delta,\theta, L) \coloneqq \sum_{v\in G_{[r]}} \mathrm e^{-\theta A_r(\delta,v, L)}.
\end{align*}
Then, the Markov inequality implies
\begin{align}\label{bound-probB^c}\notag
\p(B_r(\delta,\beta, L)^c|\xi(\rho,\rho^*)=1)
&= \p(\exists v\in G_{[r]}: A_r(\delta,v, L) < \beta r)\notag\\
&= \p\big(\min_{v\in G_{[r]}} A_r(\delta,v, L) < \beta r \big)\notag\\
&\le  \e^{\theta \beta r} \E\big[\e^{-\theta \min_{v} A_r(\delta,v, L)} \big]\notag\\
&\le \e^{\theta \beta r} \E[Z_r(\delta,\theta, L)].
\end{align}
For the expectation of $Z_r(\delta,\theta,L)$ we obtain the following recursion
\begin{align*}
\E[Z_r(\delta,\theta, L)] 
&= \E\Bigg[ \sum_{v\in G_{[r-1]}} \sum_{z\in G_L(T(v))}  \mathrm e^{-\theta A_r(\delta,z, L)}\Bigg]\\
&= \E\Bigg[ \sum_{v\in G_{[r-1]}} \mathrm e^{-\theta A_{r-1}(\delta,v, L)} \sum_{z\in G_L(T(v))} \mathrm e^{-\theta\mathds 1_{\{P^{z}_{\omega\backslash \omega(z)}(\eta_{v}=\infty)\ge \delta\}}} \Bigg]\\
&= \E\Bigg[ \sum_{v\in G_{[r-1]}} \E\Bigg[\mathrm e^{-\theta A_{r-1}(\delta,v, L)} \sum_{z\in G_L(T(v))} \mathrm e^{-\theta\mathds 1_{\{P^{z}_{\omega\backslash \omega(z)}(\eta_{v}=\infty)\ge \delta\}}}\Bigg| G_{[r-1]} \Bigg]\Bigg]\\
&= \E\Bigg[ \sum_{v\in G_{[r-1]}} \E\Bigg[\mathrm e^{-\theta A_{r-1}(\delta,v, L)} \E\Bigg[\sum_{z\in G_L(T(v))} \mathrm e^{-\theta\mathds 1_{\{P^{z}_{\omega\backslash \omega(z)}(\eta_{v}=\infty)\ge \delta\}}}\Bigg| \omega\backslash\omega(v)\Bigg] \Bigg| G_{[r-1]} \Bigg]\Bigg]\\
&= \E\Bigg[ \sum_{v\in G_{[r-1]}} \E\Bigg[\mathrm e^{-\theta A_{r-1}(\delta,v, L)} \E\Bigg[\sum_{z\in G_L(T(v))} \mathrm e^{-\theta\mathds 1_{\{P^{z}_{\omega\backslash \omega(z)}(\eta_{v}=\infty)\ge \delta\}}} \Bigg] \Bigg| G_{[r-1]} \Bigg]\Bigg]\\
&= \E\Bigg[\sum_{z\in G_L} \mathrm e^{-\theta\mathds 1_{\{P^{z}_{\omega\backslash \omega(z)}(\eta_{\rho}=\infty)\ge \delta\}}} \Bigg] \E\Bigg[ \sum_{v\in G_{[r-1]}} \E\Bigg[\mathrm e^{-\theta A_{r-1}(\delta,v, L)}  \Bigg| G_{[r-1]} \Bigg]\Bigg]\\
&=\E[Z_1(\delta,\beta, L)]\E[Z_{r-1}(\delta,\theta, L)].
\end{align*}
Iterating this gives rise to
\begin{align}\label{expectation_Z_n}
\E[Z_r(\delta,\theta, L)] = \E[Z_1(\delta,\theta, L)]^r.
\end{align}
We recall that $N$ denotes the number of offspring of the root and $m=E[N]$ is the offspring mean. Then, we can bound the expectation of $Z_1(\delta,\theta, L)$ by
\begin{align*}
\E[Z_1(\delta,\theta, L)]
&= \E\Bigg[\sum_{z\in G_L} \Big(\mathrm e^{-\theta} \mathds 1_{\{P^{z}_{\omega\backslash \omega(z)}(\eta_{\rho}=\infty)\ge \delta\}} + \mathds 1_{\{P^{z}_{\omega\backslash \omega(z)}(\eta_{\rho}=\infty)< \delta\}}\Big)  \Bigg]\\
&\le \mathrm e^{-\theta} m^L + \E\Bigg[ \sum_{z\in G_L}\mathds 1_{\{P^{z}_{\omega\backslash \omega(z)}(\eta_{\rho}=\infty)< \delta\}} \Bigg]\\
&=\mathrm e^{-\theta} m^L + \E\Bigg[ \sum_{i=1}^N \E\Bigg[ \sum_{v\in G_{L-1}(T(z_i))} \mathds 1_{\{P^{v}_{\omega^*(z_i)\backslash \omega(v)}(\eta_{\rho}=\infty)< \delta\}}\Bigg| N \Bigg]\Bigg]\\
&=\mathrm e^{-\theta} m^L + m\E\Bigg[ \sum_{v\in G_{L-1}} \mathds 1_{\{P^{v}_{\omega\backslash \omega(v)}(\eta_*=\infty)< \delta\}}\Bigg],
\end{align*}
where the descendants of the root are denoted by $z_i$. 
By Lemma~\ref{lem:quenchedindirectescapeprob} there exists a constant $C=C(\nu,\alpha,\kappa)>0$ and a function $\gamma$ with $\gamma(L,\nu,\alpha)\to 0$ for $L\to\infty$, such that
\begin{align*}
\E[Z_1(\delta,\theta, L)]
\le \mathrm e^{-\theta} m^L + m\gamma(L-1,\nu,\alpha) + C\frac\delta{1-\delta(1+L)} m^{L+1} \eqqcolon \zeta(\theta, \delta, L,\kappa).
\end{align*}
Hence, we can fix an $L$, such that $m\gamma(L,\nu,\alpha)<1$ holds. Then, we can choose $\theta$ large and $\delta$ small enough, such that the right-hand side is strictly less than one. Finally, due to \eqref{bound-probB^c} and \eqref{expectation_Z_n}, we obtain
\begin{align*}
\p(B_r(\delta,\beta, L)^c|\xi(\rho,\rho^*)=1)
\le \big(\e^{\theta \beta} \zeta(\theta, \delta, L,\kappa)\big)^r.
\end{align*}
If $\beta$ is sufficiently small, this bound decays exponentially in $r$ and then also in $n$, which completes the proof. To see this let us briefly summarize the results. We recall that, due to \eqref{bound-backtrack-2}, the estimate above implies that the event $\{\tilde N_1=n,R_1<\infty\}$ has exponentially small probability, which yields
\begin{align*}
\bbE\big[\tilde N_1^{2q}\mathds 1_{\{R_1<\infty\}}\big]\le C_2<\infty
\end{align*}
(see \eqref{moments_exptildeN}).
Hence, in view of \eqref{momentbountdist-2}, from this and the uniform bound in \eqref{momentbound-H} we conclude that $\bbE\big[|X_{\tau_1}|^q\big]$ is finite.
\hfill $\Box$

\subsection{Proof of Lemma~\ref{lem:quenchedindirectescapeprob}}
A quenched bound for the indirect escape probability highly depends on the realization of the tree. In particular, escape probability depends on the structure of the backbone tree $T_1^{\Bb}$. For this reason, in order to bound the expected number of vertices with an indirect escape probability less than $\delta$, we distinguish the different sizes of the $L$-th generation of $T_1^{\Bb}$
\begin{align}\label{expectation_decompositionBbSizes}
\E\Bigg[ \sum_{v\in G_L}\mathds 1_{\{P^{v}_{\omega\backslash \omega(v)}(\eta_*=\infty)< \delta\}} \Bigg] 
&= \E\Bigg[ \sum_{v\in G_L}\mathds 1_{\{P^{v}_{\omega\backslash \omega(v)}(\eta_*=\infty)< \delta\}}\mathds 1_{\{|G_L\cap T_1^{\Bb}|=0\}} \Bigg] \notag\\
&\quad+ \E\Bigg[ \sum_{v\in G_L}\mathds 1_{\{P^{v}_{\omega\backslash \omega(v)}(\eta_*=\infty)< \delta\}}\mathds 1_{\{|G_L\cap T_1^{\Bb}|=1\}} \Bigg] \notag\\
&\quad+ \E\Bigg[ \sum_{v\in G_L}\mathds 1_{\{P^{v}_{\omega\backslash \omega(v)}(\eta_*=\infty)< \delta\}}\mathds 1_{\{|G_L\cap T_1^{\Bb}|>1\}} \Bigg]
\end{align}
We study the three summands separately.

\subsubsection{The indirect escape probability in the case that the $L$-th generation of $T_1^{\Bb}$ contains zero vertices}
We start with considering the expected number of vertices with an indirect escape probability less than $\delta$ on the event, that the $L$-th generation of the backbone tree $T_1^{\Bb}$ is empty, which means that the tree $T_1(\rho)$ is finite.

For any vertex $v\in G_L$ we denote its ancestor in the $i$-th generation $G_i$ by $v^i$, i.e. $\rho^*=v^{-1},\rho=v^0,v^1,\dots,v^L=v$ is the unique path from $\rho^*$ to $v$. If $v\notin T_1(\rho)$, then we let $k_\varepsilon$ be the index, such that $(v^{k_\varepsilon-1},v^{k_\varepsilon})$ is the first edge on the path from $\rho^*$ to $v$ with conductance $\varepsilon$, i.e. 
\begin{align*}
k_\varepsilon 
\coloneqq k_\varepsilon(v) 
\coloneqq \min\{ k: \xi(v^{k-1},v^k)=\varepsilon\}.
\end{align*}
Further, we introduce the event
\begin{align*}
A \coloneqq \{v\in T\backslash T_1(\rho) : \deg(v^k)=2 \text{ for all } k\in\{k_\varepsilon(v),\dots,|v|-1\}.
\end{align*}
\begin{figure}
	\centering
	\begin{tikzpicture}
	\draw (-1.7,0) node[] {$\rho^*$};
	\draw[-, line width=0.5mm] (-1.5,0) to (0,0);
	\draw (-0.1,-0.2) node[] {$\rho$}; 
	\draw[-, line width=0.5mm] (0,0) to (1.5,1.5);
	\draw[-, line width=0.5mm] (0,0)to (1.5,-1.5);
	\draw[-, line width=0.5mm] (1.5,1.5) to (3,2.5);
	\draw[style=dotted] (1.5,-1.5) to (3,-0.7);
	\draw[-, line width=0.5mm] (1.5,-1.5) to (3,-2.5);
	\draw[-, line width=0.5mm] (3,2.5) to (4.5,3.25); 
	\draw[-, line width=0.5mm] (3,2.5) to (4.5,1.75); 
	\draw[style=dotted] (3,-0.7) to (4.5,0); 
	\draw[style=dotted] (3,-0.7) to (4.5,-1.4); 
	\draw[style=dotted] (3,-2.5) to (4.5,-3.15); \draw (4.5,-3.55) node[] {$v_2^{k_\varepsilon}$}; \draw[] (4.5,-3.15) circle(3pt);
	\draw[style=dotted] (4.5,3.25) to (6,3.75); \draw (6,4.15) node[] {$v_1^{k_\varepsilon}$}; \draw[] (6,3.75) circle(3pt);
	\draw[style=dotted] (4.5,3.25) to (6,2.75);
	\draw[-, line width=0.5mm] (4.5,1.75) to (6,1.75);  
	\draw[-] (4.5,0) to (6,0); 
	\draw[-] (4.5,-1.4) to (6,-0.9); 
	\draw[style=dotted] (4.5,-1.4) to (6,-1.9); 
	\draw[style=dotted] (4.5,-3.15) to (6,-3.5); \draw[] (6,-3.5) circle(3pt);
	\draw[-] (6,3.75) to (7.5,4);  \draw[] (7.5,4) circle(3pt);
	\draw[style=dotted] (6,2.75) to (7.5,3.05);
	\draw[style=dotted] (6,2.75) to (7.5,2.45);
	\draw[-, line width=0.5mm] (6,1.75) to (7.5,2.05); 
	\draw[-, line width=0.5mm] (6,1.75) to (7.5,1.45); 
	\draw[style=dotted] (6,0) to (7.5,0.3); 
	\draw[-] (6,0) to (7.5,-0.3);
	\draw[-] (6,-0.9) to (7.5,-0.9); 
	\draw[-] (6,-1.9) to (7.5,-1.6); 
	\draw[-] (6,-1.9) to (7.5,-2.2); 
	\draw[-] (6,-3.5) to (7.5,-3.75); \draw[] (7.5,-3.75) circle(3pt);
	\draw[-] (7.5,4) to (9,4); \fill (9,4) circle (2pt); \draw (9,4) node[right] {\small $v_1\in A$};
	\draw[style=dotted] (7.5,3.05) to (9,3.25); \fill (9,3.25) circle (2pt); 
	\draw[style=dotted] (7.5,3.05) to (9,2.85); \fill (9,2.85) circle (2pt); 
	\draw[-] (7.5,2.45) to (9,2.45); \fill (9,2.45) circle (2pt); 
	\draw[-, line width=0.5mm] (7.5,2.05) to (9,2.05); \fill (9,2.05) circle (2pt); \draw (9,2.05) node[right] {\small $z_1\in T_1(\rho)$};
	\draw[-, line width=0.5mm] (7.5,1.45) to (9,1.65); \fill (9,1.65) circle (2pt); \draw (9,1.65) node[right] {\small $z_2\in T_1(\rho)$};
	\draw[-, line width=0.5mm] (7.5,1.45) to (9,1.25); \fill (9,1.25) circle (2pt); \draw (9,1.25) node[right] {\small $z_3\in T_1(\rho)$};
	\draw[-] (7.5,0.3) to (9,0.3); \fill (9,0.3) circle (2pt); 
	\draw[-] (7.5,-0.3) to (9,-0.1); \fill (9,-0.1) circle (2pt); 
	\draw[-] (7.5,-0.3) to (9,-0.5); \fill (9,-0.5) circle (2pt); 
	\draw[-] (7.5,-0.9) to (9,-0.9); \fill (9,-0.9) circle (2pt); 
	\draw[-] (7.5,-1.6) to (9,-1.6); \fill (9,-1.6) circle (2pt); 
	\draw[style=dotted] (7.5,-2.2) to (9,-2); \fill (9,-2) circle (2pt); 
	\draw[-] (7.5,-2.2) to (9,-2.4); \fill (9,-2.4) circle (2pt); 
	\draw[-] (7.5,-3.75) to (9,-3.75); \fill (9,-3.75) circle (2pt); \draw (9,-3.75) node[right] {\small $v_{2}\in A$};
	\draw[dashed] (9,4.5) to (9,-4.5);
	\draw (9.65,-4.5) node[] {$G_6(T)$};
	\end{tikzpicture}
	\caption{ 
		Edges with conductance larger $\varepsilon$ are indicated by solid lines; edges with conductance $\varepsilon$ are indicated by dotted lines; vertices in the sixth generation $G_6(T)$ are marked by dots on the dashed line.
		$T_1(\rho)$ is the subtree formed by the edges with conductance one containing the root (indicated by the thick lines). In the sixth generation $G_6(T)$ the vertices $z_1,\,z_2,\,z_3$ are in $T_1(\rho)$.
		The event $A$ contains a vertex $v\in T$ if every vertex on the path from $v^{k_\varepsilon}$ to $v^*$ has degree 2, where $(v^{k_\varepsilon-1},v^{k_\varepsilon})$ is the first $\varepsilon$-edge on the path from $\rho^*$ to $v$. In the sixth generation $G_6(T)$ the vertices $v_1,\,v_2$ are in $A$. The ancestors of $v_i$, which must have degree 2 for $v_i$ being in $A$, are circled.	
	}
	\label{fig:eventA}
\end{figure}
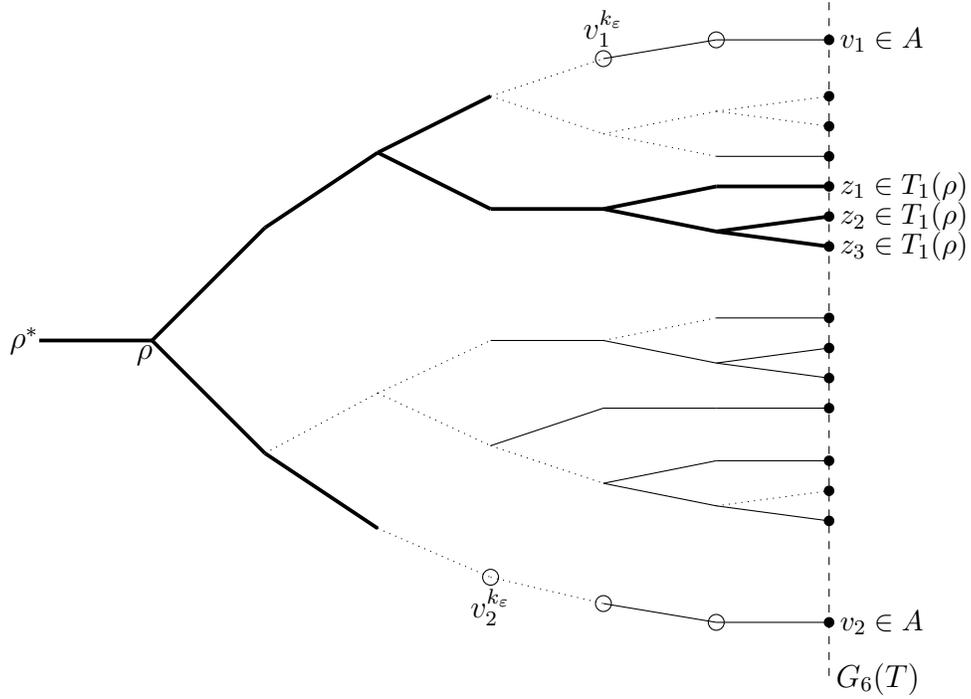

In order to compute the first expectation of \eqref{expectation_decompositionBbSizes}, we distinguish whether the vertex $v$ is in $T_1(\rho)$, in $A$ or in neither (see figure~\ref{fig:eventA} for an example)
\begin{align} \label{backbone=0}
\E\Bigg[ \sum_{v\in G_L}\mathds 1_{\{P^{v}_{\omega\backslash \omega(v)}(\eta_*=\infty)< \delta\}}\mathds 1_{\{|G_L\cap T_1^{\Bb}|=0\}} \Bigg] 
&= \E\Bigg[ \sum_{v\in G_L}\mathds 1_{\{P^{v}_{\omega\backslash \omega(v)}(\eta_*=\infty)< \delta\}}\mathds 1_{\{|T_1(\rho)|<\infty\}} \mathds 1_{\{v\in T_1(\rho)\}} \Bigg] \notag \\
&+  \E\Bigg[ \sum_{v\in G_L}\mathds 1_{\{P^{v}_{\omega\backslash \omega(v)}(\eta_*=\infty)< \delta\}}\mathds 1_{\{|T_1(\rho)|<\infty\}} \mathds 1_{\{v\in A\}} \Bigg] \notag \\
&+ \E\Bigg[ \sum_{v\in G_L}\mathds 1_{\{P^{v}_{\omega\backslash \omega(v)}(\eta_*=\infty)< \delta\}}\mathds 1_{\{|T_1(\rho)|<\infty\}} \mathds 1_{\{v\notin T_1(\rho)\cup A\}} \Bigg].
\end{align}

We will treat the three summands separately in order of their appearance. \\

\textbf{Bound for first expectation of \eqref{backbone=0}:}
When the random walk starts in $v\in T_1(\rho)$, we can not expect a small probability to never reach $\rho^*$, because in this situation the walker has to take an $\varepsilon$-edge to escape.
Therefore, we bound the first expectation of \eqref{backbone=0} by the expected size of the $L$-th generation of the subtree $T_1(\rho)$ conditioned on this tree to die out. That is
\begin{align}\label{backbone=0_boundEW1}
\E\Bigg[ \sum_{v\in G_L}\mathds 1_{\{P^{v}_{\omega\backslash \omega(v)}(\eta_*=\infty) < \delta\}}\mathds 1_{\{|T_1(\rho)|<\infty\}} \mathds 1_{\{v\in T_1(\rho)\}} \Bigg]
&\le \E\Bigg[ \sum_{v\in G_L} \mathds 1_{\{v\in T_1(\rho)\}} \Bigg| |T_1(\rho)|<\infty \Bigg] \notag\\
&= (m_1')^L
\end{align}
with $m_1'\coloneqq \E\big[\sum_{z\in G_1} \mathds 1_{\{\kappa^{-1}\le\xi(\rho,z)\le\kappa\}} \big| |T_1(\rho)|<\infty \big]$. The subtree $T_1(\rho)$, conditioned on extinction, is a subcritical Galton-Watson tree, which implies that $m_1'$ is strictly less than one (see \cite{hofstad2016randomgraphs} Theorem~3.7 and Exercise~3.17). This shows that the first expectation on the right-hand side of \eqref{backbone=0} vanishes as $L$ tends to infinity. \\

\textbf{Bound for second expectation of \eqref{backbone=0}:}
We can bound the second expectation in \eqref{backbone=0} by the expected number of vertices $v\in G_L$ in $A$
\begin{align*}
\E\Bigg[ \sum_{v\in G_L}\mathds 1_{\{P^{v}_{\omega\backslash \omega(v)}(\eta_*=\infty)< \delta\}}\mathds 1_{\{|T_1(\rho)|<\infty\}} \mathds 1_{\{v\in A\}} \Bigg]
\le \E\Bigg[ \sum_{v\in G_L} \mathds 1_{\{|T_1(\rho)|<\infty\}} \mathds 1_{\{v\in A\}} \Bigg].
\end{align*}
In order to compute this expectation, we distinguish whether the ancestor $v^*$ of $v$ is in $A$ as well or in $T_1(\rho)$
\begin{align}\label{backbone=0_EW2}
&\E\Bigg[ \sum_{v\in G_L} \mathds 1_{\{|T_1(\rho)|<\infty\}} \mathds 1_{\{v\in A\}} \Bigg] 
= \E\Bigg[ \sum_{v\in G_L} \mathds 1_{\{|T_1(\rho)|<\infty\}} \big(\mathds 1_{\{v^*\in A,\,\deg(v^*)=2\}}+ \mathds 1_{\{v^*\in T_1(\rho),\,\xi(v^*,v)=\varepsilon\}}\big) \Bigg] \notag\\
&=\E\Bigg[ \sum_{v\in G_{L-1}} \mathds 1_{\{|T_1(\rho)|<\infty\}}\mathds 1_{\{v\in A,\,\deg(v)=2\}} \Bigg]
+ \E\Bigg[ \sum_{v\in G_{L-1}} \mathds 1_{\{|T_1(\rho)|<\infty\}} \sum_{z\in G_1(T(v))} \mathds 1_{\{v\in T_1(\rho),\,\xi(v,z)=\varepsilon\}} \Bigg].
\end{align}
The first expectation in \eqref{backbone=0_EW2} is given by
\begin{align}\label{backbone=0_EW2_A}
\E\Bigg[ \sum_{v\in G_{L-1}} \mathds 1_{\{|T_1(\rho)|<\infty\}}\mathds 1_{\{v\in A,\,\deg(v)=2\}} \Bigg] 
&= \E\Bigg[ \sum_{v\in G_{L-1},\,v\in A} \E\big[ \mathds 1_{\{|T_1(\rho)|<\infty\}} \mathds 1_{\{\deg(v)=2\}} \big| \omega_{|L-1} \big]\Bigg] \notag\\
&= \E\Bigg[\sum_{v\in G_{L-1},\,v\in A} \p(\deg(v)=2) \E\big[ \mathds 1_{\{|T_1(\rho)|<\infty\}} \big| \omega_{|L-1} \big] \Bigg]\notag\\
&= p_1 \E\Bigg[\sum_{v\in G_{L-1}} \mathds 1_{\{v\in A\}}\mathds 1_{\{|T_1(\rho)|<\infty\}} \Bigg],
\end{align}
where $\omega_{|L-1}$ denotes the first $L-1$ generations of the environment and $p_1=\nu(\{1\})$ is the probability that a vertex has one descendant.
To see the second equality, observe that $\mathds 1_{\{\deg(v)=2\}}$ and $\mathds 1_{\{|T_1(\rho)|<\infty\}}$ are independent, because the first indicator function is $\sigma(\omega(v))$-measurable and the second one is $\sigma(\omega\backslash\omega(v))$-measurable, since $v\notin T_1(\rho)$.
For the second summand in \eqref{backbone=0_EW2} we obtain
\begin{align*}
&\E\Bigg[ \sum_{v\in G_{L-1}} \mathds 1_{\{|T_1(\rho)|<\infty\}} \sum_{z\in G_1(T(v))} \mathds 1_{\{v\in T_1(\rho),\,\xi(v,z)=\varepsilon\}} \Bigg]\\
&\le \E\Bigg[ \sum_{\substack{v\in G_{L-1}, v\in T_1(\rho)}} \E\bigg[ \mathds 1_{\{|T_1(y)|<\infty \text{ for all } y\in G_{L-1}\cap T_1(\rho)\setminus\{v\} \}} \sum_{z\in G_1(T(v))} \mathds 1_{\{\xi(v,z)=\varepsilon\}}\bigg| \omega_{|L-1}\bigg] \Bigg].
\end{align*}
Here we used that the tree $T_1(\rho)$ can only die out when the subtrees $T_1(y)$ become extinct for all $y\in G_{L-1}$ in the $(L-1)$-th generation of $T_1(\rho)$. Hence, we get an upper bound if this holds for all $y\in G_{L-1}\cap T_1(\rho)$ except $v$.
The indicator function $\mathds 1_{\{|T_1(y)|<\infty \text{ for all } y\in G_{L-1}\cap T_1(\rho)\setminus\{v\} \}}$ is $\sigma(\omega\backslash\omega(v))$-measurable whereas $\mathds 1_{\{|T_1(v)|<\infty\}}$ and \linebreak $\sum_{z\in G_1(T(v))} \mathds 1_{\{\xi(v,z)=\varepsilon\}}$ are $\sigma(\omega(v))$-measurable. Therefore, they are independent, which implies
\begin{align*}
&\E\Bigg[ \sum_{v\in G_{L-1}} \mathds 1_{\{|T_1(\rho)|<\infty\}} \sum_{z\in G_1(T(v))} \mathds 1_{\{\xi(v,z)=\varepsilon\}} \Bigg]\\
&=\E\Bigg[ \sum_{\substack{v\in G_{L-1}, v\in T_1(\rho)}} \E\big[ \mathds 1_{\{|T_1(y)|<\infty \text{ for all } y\in G_{L-1}\cap T_1(\rho)\setminus\{v\} \}}\big| \omega_{|L-1} \big] 
\E\bigg[ \sum_{z\in G_1(T(v))} \mathds 1_{\{\xi(v,z)=\varepsilon\}} \bigg]\Bigg]\\
&=m_\varepsilon \E\Bigg[ \sum_{\substack{v\in G_{L-1}, v\in T_1(\rho)}}  \E\big[ \mathds 1_{\{|T_1(y)|<\infty \text{ for all } y\in G_{L-1}\cap T_1(\rho)\setminus\{v\} \}}\big| \omega_{|L-1}\big]\Bigg]\\
&=m_\varepsilon\E\Bigg[ \sum_{\substack{v\in G_{L-1}, v\in T_1(\rho)}}  \E\big[ \mathds 1_{\{|T_1(y)|<\infty \text{ for all } y\in G_{L-1}\cap T_1(\rho) \}}\big| \omega_{|L-1}\big]\p(|T_1(v)|<\infty)^{-1}\Bigg]\\
&=m_\varepsilon \E\Bigg[ \sum_{v\in G_{L-1}}  \mathds 1_{\{v\in T_1(\rho)\}} \mathds 1_{\{|T_1(\rho)|<\infty \}} \Bigg] \p(|T_1(\rho)|<\infty)^{-1}\\
&=m_\varepsilon \E\Bigg[ \sum_{v\in G_{L-1}}  \mathds 1_{\{v\in T_1(\rho)\}} \Bigg| |T_1(\rho)|<\infty \Bigg],
\end{align*}
where $m_\varepsilon = \E\big[\sum_{z\in G_1} \mathds 1_{\{\xi(\rho,z)=\varepsilon\}}\big]$ denotes the expected number of descendants of the root which are connected by an $\varepsilon$-edge. The subtree $T_1(\rho)$,  conditioned on extinction, is a subcritical Galton-Watson tree. Hence, conditioned on the extinction of $T_1(\rho)$, the expected size of the $(L-1)$-th generation of $T_1(\rho)$ is given by $(m_1')^{L-1}$. Consequently, we have
\begin{align}\label{backbone=0_EW2_B}
\E\Bigg[ \sum_{v\in G_{L-1}} \mathds 1_{\{|T_1(\rho)|<\infty\}} \sum_{z\in G_1(T(v))} \mathds 1_{\{\xi(v,z)=\varepsilon\}} \Bigg]
\le m_\varepsilon (m_1')^{L-1}.
\end{align} 
Plugging \eqref{backbone=0_EW2_A} and \eqref{backbone=0_EW2_B} in \eqref{backbone=0_EW2}, we obtain 
\begin{align}\label{backbone=0_EW2_recursion}
\E\Bigg[ \sum_{v\in G_L} \mathds 1_{\{|T_1(\rho)|<\infty\}} \mathds 1_{\{v\in A\}} \Bigg] 
\le p_1 \E\Bigg[\sum_{v\in G_{L-1}} \mathds 1_{\{v\in A\}}\mathds 1_{\{|T_1(\rho)|<\infty\}} \Bigg] 
+ m_\varepsilon (m_1')^{L-1}.
\end{align}
Hence, we have a recursion of the form $x_{n+1}\le ax_n+cb^n$, $x_0=0$, for some $a,\,b,\,c\ge 0$. Iterating this gives rise to $x_{n+1}\le c\sum_{i=0}^n a^ib^{n-i}$. In our setting, iterating \eqref{backbone=0_EW2_recursion} leads to
\begin{align}\label{backbone=0_boundEW2}
\E\Bigg[ \sum_{v\in G_L} \mathds 1_{\{|T_1(\rho)|<\infty\}} \mathds 1_{\{v\in A\}} \Bigg] 
\le m_\varepsilon \sum_{i=0}^{L-1} p_1^i(m_1')^{L-1-i}
= m_\varepsilon (m_1')^{L-1} \sum_{i=0}^{L-1} \Big(\frac{p_1}{m_1'}\Big)^{i}.
\end{align}
which vanishes as $L\to \infty$, since $m_1',p_1<1$.  \\

\textbf{Bound for third expectation of \eqref{backbone=0}:}
Roughly speaking, when the random walk starts in a vertex $v\notin T_1(\rho)\cup A$, it has a good probability to never visit $\rho^*$. Because on the one hand starting in a vertex outside of $T_1(\rho)$ implies that the walker has to move along an $\varepsilon$-edge to reach $\rho^*$ at some time, one the other hand there exists a path in $T\backslash T(v)$ to escape, since the random walk starts in a vertex that is not located in the set $A$. For this reason, to bound the third expectation of \eqref{backbone=0}, we deduce a lower bound for the indirect escape probability.

We let $\omega$ be an environment with $|G_L\cap T_1^{\Bb}|=0$ and $v\in G_L$ a vertex in the $L$-th generation with $v\notin T_1(\rho)$ and $v\notin A$. Recall that $(v^{k_\varepsilon-1},v^{k_\varepsilon})$ denotes the first $\varepsilon$-edge on the path from $\rho^*$ to $v$, we have
\begin{align*}
P^v_{\omega\backslash\omega(v)}(\eta_*=\infty)
\ge P^{v^{k_\varepsilon}}_{\omega\backslash\omega(v)} (\eta_{v^{k_\varepsilon-1}}=\infty)
=\frac{\calC_{\omega(v^{k_\varepsilon})\backslash \omega(v)}(v^{k_\varepsilon},\infty)}{\varepsilon + \calC_{\omega(v^{k_\varepsilon})\backslash \omega(v)}(v^{k_\varepsilon},\infty)}.
\end{align*}
By Rayleigh's monotonicity principle, the effective conductance decreases when the conductances of all edges are reduced to $\varepsilon$, which implies
\begin{align*}
\calC_{\omega\backslash \omega(v)}(v^{k_\varepsilon},\infty) \ge \calC_{\omega_{[\varepsilon]}(v^{k_\varepsilon})\backslash \omega_{[\varepsilon]}(v)}(v^{k_\varepsilon},\infty).
\end{align*}
Here $\omega_{[a]}=(T,\rho,(a)_{e\in\calE(T)})$ denotes the environment, where each conductance in $\omega$ is replaced by $a>0$. The monotonicity of the mapping $x\mapsto \frac x{x+\varepsilon}$ then implies
\begin{align*}
P^v_{\omega\backslash\omega(v)}(\eta_*=\infty)
&\ge \frac{\calC_{\omega_{[\varepsilon]}(v^{k_\varepsilon})\backslash \omega_{[\varepsilon]}(v)}(v^{k_\varepsilon},\infty)}{\varepsilon + \calC_{\omega_{[\varepsilon]}(v^{k_\varepsilon})\backslash \omega_{[\varepsilon]}(v)}(v^{k_\varepsilon},\infty)}
= P^{v^{k_\varepsilon}}_{\omega_{[\varepsilon]}\backslash\omega_{[\varepsilon]}(v)} (\eta_{v^{k_\varepsilon-1}}=\infty)\\
&= P^{v^{k_\varepsilon}}_{\omega_{[1]}\backslash\omega_{[1]}(v)} (\eta_{v^{k_\varepsilon-1}}=\infty) 
=\frac{\calC_{\omega_{[1]}(v^{k_\varepsilon})\backslash \omega_{[1]}(v)}(v^{k_\varepsilon},\infty)}{1 + \calC_{\omega_{[1]}(v^{k_\varepsilon})\backslash \omega_{[1]}(v)}(v^{k_\varepsilon},\infty)}\\
&=\frac 1{1+\calR_{\omega_{[1]}\backslash \omega_{[1]}(v)}(v^{k_\varepsilon},\infty)}.
\end{align*}
Since $v\notin A$, there exists an index $k'\in\{k_\varepsilon,\dots,L-1\}$, such that $\deg(v^{k'})>2$. In particular, there exists a vertex $u\in G_L\cap T(v^{k_\varepsilon})\backslash T(v^{k'+1})$.
Again by Rayleigh's monotonicity principle the effective resistance can only increase, when we remove edges. Together with the Series Law we get
\begin{align*}
\calR_{\omega_{[1]}\backslash \omega_{[1]}(v)}(v^{k_\varepsilon},\infty) 
\le \calR_{\omega_{[1]}}(v^{k_\varepsilon},u) + \calR_{\omega_{[1]}(u)}(u,\infty)
\le \kappa L + \calR_{\omega_{[1]}(u)}(u,\infty),
\end{align*}
thus
\begin{align*}
P^v_{\omega\backslash\omega(v)}(\eta_*=\infty) 
\ge \frac 1{1+\kappa L + \calR_{\omega_{[1]}(u)}(u,\infty)}.
\end{align*}
From this we conclude that $ P^v_{\omega\backslash\omega(v)}(\eta_*=\infty) < \delta$ can only hold, if 
\begin{align*}
\calR_{\omega_{[1]}(u)}(u,\infty) > \frac {1-\delta(1+\kappa L)}\delta.
\end{align*}
This allows us to estimate the third expectation of \eqref{backbone=0} as follows
\begin{align}\label{backbone=0_EW3}
&\E\Bigg[ \sum_{v\in G_L}\mathds 1_{\{P^{v}_{\omega\backslash \omega(v)}(\eta_*=\infty)< \delta\}}\mathds 1_{\{|T_1(\rho)|<\infty\}} \mathds 1_{\{v\notin T_1(\rho)\cup A\}} \Bigg]\notag\\
&\quad\le \E\Bigg[ \sum_{\substack{v\in G_{L}}} \mathds 1_{\{v\notin T_1(\rho)\cup A\}} \E\big[\mathds 1_{\{P^{v}_{\omega\backslash \omega(v)} (\eta_*=\infty)< \delta\}}\big|G_{L-1} \big]\Bigg]\notag\\
&\quad\le \E\Bigg[ \sum_{\substack{v\in G_{L}}} \mathds 1_{\{v\notin T_1(\rho)\cup A\}} \E\bigg[ \mathds 1_{\big\{\calR_{\omega_{[1]}(u)}(u,\infty) > \frac {1-\delta(1+\kappa L)}\delta\big\}} \bigg]\Bigg]\notag\\
&\quad= \p\bigg(\calR_{\omega_{[1]}}(\rho,\infty) > \frac {1-\delta(1+\kappa L)}\delta\bigg) \E\Bigg[ \sum_{v\in G_L} \mathds 1_{\{v\notin T_1(\rho)\cup A\}} \Bigg].
\end{align}
For $\delta>0$ small enough, applying the Markov inequality gives rise to the following bound
\begin{align}\label{backbone=0_EW3_A}
\p\bigg(\calR_{\omega_{[1]}}(\rho,\infty) > \frac {1-\delta(1+\kappa L)}\delta\bigg)
\le \frac \delta{1-\delta(1+\kappa L)} \E\big[\calR_{\omega_{[1]}}(\rho,\infty)\big].
\end{align}
Since $\omega_{[1]}$ is a supercritical Galton-Watson tree with unit conductance, from \cite[Lemma~9.1]{LyoPemPer96a} follows, that the expected effective resistance is finite. Further, we can bound the expectation in \eqref{backbone=0_EW3} by the expected size of the $L$-th generation
\begin{align*}
\E\Bigg[ \sum_{v\in G_L} \mathds 1_{\{v\notin T_1(\rho)\cup A\}} \Bigg] 
\le E[|G_L|]=m^L,
\end{align*}
recall that $m$ is the offspring mean.  Plugging this estimate and \eqref{backbone=0_EW3_A} in \eqref{backbone=0_EW3}, we get
\begin{align}\label{backbone=0_boundEW3}
\E\Bigg[ \sum_{v\in G_L}\mathds 1_{\{P^{v}_{\omega\backslash \omega(v)}(\eta_*=\infty)< \delta\}}\mathds 1_{\{|T_1(\rho)|<\infty\}} \mathds 1_{\{v\notin T_1(\rho)\cup A\}} \Bigg]
\le C_1 \frac\delta{1-\delta(1+\kappa L)} m^L
\end{align}
for some constant $C_1=C_1(\nu)>0$. 

Finally, due to \eqref{backbone=0_boundEW1}, \eqref{backbone=0_boundEW2} and \eqref{backbone=0_boundEW3}, we obtain the following bound for the expectation in \eqref{backbone=0}
\begin{align}\label{backbone=0_bound}
\E\Bigg[ \sum_{v\in G_L}\mathds 1_{\{P^{v}_{\omega\backslash \omega(v)}(\eta_*=\infty)< \delta\}}\mathds 1_{\{|G_L\cap T_1^{\Bb}|=0\}} \Bigg] 
\le h(L,\nu,\alpha) 
+ C_1\frac\delta{1-\delta(1+\kappa L)} m^L
\end{align}
for some function $h$ that vanishes as $L$ tends to infinity. This bound is independent of $\varepsilon$ and becomes arbitrarily small, when we first choose $L$ large enough and afterwards $\delta$ sufficiently small. 

We remark that for a subcritical tree $T_1(\rho)$, we have $G_L\cap T_1^{\Bb}=\emptyset$ and therefore
\begin{align*}
	\E\Bigg[ \sum_{v\in G_L}\mathds 1_{\{P^{v}_{\omega\backslash \omega(v)}(\eta_*=\infty)< \delta\}}\mathds 1_{\{|G_L\cap T_1^{\Bb}|\ge 1\}} \Bigg] =0,
\end{align*}
which shows that in the (sub-)critical case we are done at this point.
Therefore, for the rest of the proof we assume the tree $T_1(\rho)$ to be supercritical.\\

\subsubsection{The indirect escape probability in the case that the $L$-th generation of $T_1^{\Bb}$ contains one vertex}
We proceed by considering the expected number of vertices with an indirect escape probability less than $\delta$ on the event, that the $L$-th generation of $T_1^{\Bb}$ contains exactly one vertex. We distinguish whether the random walk starts at this vertex or not
\begin{align} \label{backbone=1}
\E\Bigg[ \sum_{v\in G_L}\mathds 1_{\{P^{v}_{\omega\backslash \omega(v)}(\eta_*=\infty)< \delta\}}\mathds 1_{\{|G_L\cap T_1^{\Bb}|=1\}} \Bigg] 
&= \E\Bigg[ \sum_{v\in G_L}\mathds 1_{\{P^{v}_{\omega\backslash \omega(v)}(\eta_*=\infty)< \delta\}}\mathds 1_{\{|G_L\cap T_1^{\Bb}|=1\}} \mathds 1_{\{v\in T_1^{\Bb}\}} \Bigg] \notag \\
&+ \E\Bigg[ \sum_{v\in G_L}\mathds 1_{\{P^{v}_{\omega\backslash \omega(v)}(\eta_*=\infty)< \delta\}}\mathds 1_{\{|G_L\cap T_1^{\Bb}|=1\}} \mathds 1_{\{v\notin T_1^{\Bb}\}} \Bigg].
\end{align}

\textbf{Bound for the first expectation of \eqref{backbone=1}:}
The first expectation of \eqref{backbone=1} is bounded by the probability that the $L$-th generation of the backbone tree $T_1^{\Bb}$ only contains one vertex
\begin{align}\label{backbone=1_EW1}
&\E\Bigg[ \sum_{v\in G_L}\mathds 1_{\{P^{v}_{\omega\backslash \omega(v)}(\eta_*=\infty)< \delta\}}\mathds 1_{\{|G_L\cap T_1^{\Bb}|=1\}} \mathds 1_{\{v\in T_1^{\Bb}\}} \Bigg] \notag\\
&\quad\le \E\Bigg[ \sum_{v\in G_L} \mathds 1_{\{v\in T_1^{\Bb}\}} \Bigg| |G_L\cap T_1^{\Bb}|=1 \Bigg] \p\big(|G_L\cap T_1^{\Bb}|=1\big) \notag\\
&= \p\big(|G_L\cap T_1^{\Bb}|=1\big)
\le \p\big(|G_L\cap T_1^{\Bb}|=1\big| |T_1(\rho)|=\infty \big).
\end{align}
Conditioned on the survival of $T_1(\rho)$, the backbone tree $T_1^{\Bb}$ is again a supercritical Galton-Watson tree, which implies
\begin{align*}
\lim_{L\to\infty} \p\big(|G_L\cap T_1^{\Bb}|=1\big| |T_1(\rho)|=\infty \big) = 0 . 
\end{align*}

\textbf{Bound for the second expectation of \eqref{backbone=1}:}
When the random walk does not start at the vertex in $T_1^{\Bb}$, it can escape via the path that belongs to the backbone tree. This implies that the random walk has a good probability of never hitting the vertex $\rho^*$.
Just as for the third expectation in the previous case, we will derive a lower bound for the indirect escape probability to bound the second expectation of \eqref{backbone=1}.

We let $\omega$ be an environment with $|G_L\cap T_1^{\Bb}|=1$ and $v\in G_L$ a vertex in the $L$-th generation with $v\notin T_1^{\Bb}$. Then, there exists exactly one vertex $u\in G_L$ such that the subtree $T_1(u)$ survives. We can express the indirect escape probability as a ratio of effective conductances
\begin{align*}
P^{v}_{\omega\backslash \omega(v)}(\eta_*=\infty)
\ge \p^\rho_{\omega\backslash \omega(v)}(\eta_*=\infty)
=\frac{\calC_{\omega\backslash \omega(v)}(\rho,\infty)}{\xi(\rho,\rho^*)+\calC_{\omega\backslash \omega(v)}(\rho,\infty)}
\ge\frac1{1+\kappa \calR_{\omega\backslash\omega(v)}(\rho,\infty)}.
\end{align*}
By Rayleigh's monotonicity principle, the effective resistance can only increase when edges are removed. Together with the Series Law this implies
\begin{align*}
\calR_{\omega\backslash\omega(v)}(\rho,\infty)
\le \calR_\omega(\rho,u)+\calR_{\omega(u)}(u,\infty)
\le \kappa L + \calR_{\omega_1^{\Bb}(u)}(u,\infty).
\end{align*}
We have then
\begin{align}\label{backbone=1_EW2_indirectescapeprob}
P^{v}_{\omega\backslash \omega(v)}(\eta_*=\infty)
\ge \frac1{1+\kappa^2 L + \kappa \calR_{\omega_1^{\Bb}(u)}(u,\infty)},
\end{align}
which implies that $\p_{\omega\backslash\omega(v)}^v(\eta_*=\infty)<\delta$ can only be satisfied if
\begin{align*}
\calR_{\omega_1^{\Bb}(u)}(u,\infty) > \frac{1-\delta(1+\kappa^2 L)}{\delta\kappa}.
\end{align*}
We obtain
\begin{align*}
&\E\Bigg[ \sum_{v\in G_L}\mathds 1_{\{P^{v}_{\omega\backslash \omega(v)}(\eta_*=\infty)< \delta\}}\mathds 1_{\{|G_L\cap T_1^{\Bb}|=1\}} \mathds 1_{\{v\notin T_1^{\Bb}\}} \Bigg]\\
&=\E\Bigg[ \sum_{v\in G_L} \sum_{u\in G_L\cap T_1(\rho)}
\E\big[ \mathds 1_{\{P^{v}_{\omega\backslash \omega(v)}(\eta_*=\infty)< \delta\}}\mathds 1_{\{|G_L\cap T_1^{\Bb}|=1\}} \mathds 1_{\{v\notin T_1^{\Bb}\}} \mathds 1_{\{|T_1(u)|=\infty\}}\big| G_L \big]\Bigg]\\
&\le \E\Bigg[ \sum_{v\in G_L} \sum_{u\in G_L\cap T_1(\rho)}
\E\bigg[ 
\mathds 1_{\big\{ \calR_{\omega_1^{\Bb}(u)}(u,\infty) > \frac{1-\delta(1+\kappa^2L)}{\delta\kappa},\,|T_1(u)|=\infty \big\}}
\mathds 1_{\{|T_1(z)|<\infty \text{ for all }z\in G_L\cap T_1(\rho)\backslash\{u\},\,v\notin T_1^{\Bb} \} } \bigg|G_L\bigg]\Bigg]\\
&= \E\Bigg[ \sum_{v\in G_L} \sum_{u\in G_L\cap T_1(\rho)}
\E\big[ \mathds 1_{\{|T_1(z)|<\infty \text{ for all }z\in G_L\cap T_1(\rho)\backslash\{u\} ,\,v\notin T_1^{\Bb} \} } \big| G_L \big]\Bigg]
\\&\hspace{0.8cm}\times
\E\bigg[ \mathds 1_{\big\{ \calR_{\omega_1^{\Bb}(\rho)}(\rho,\infty) > \frac{1-\delta(1+\kappa^2L)} {\delta\kappa},\,|T_1(\rho)|=\infty \big\}} \bigg],
\end{align*}
where we used the independence of the indicator functions in the last equality. For $\delta>0$ sufficiently small the Markov inequality gives rise to 
\begin{align*}
\E\bigg[ \mathds 1_{\big\{ \calR_{\omega_1^{\Bb}(\rho)}(\rho\infty) > \frac{1-\delta(1+\kappa^2L)}{ \delta\kappa},\,|T_1(\rho)|=\infty \big\}} \bigg]
&\le\p\bigg(\calR_{\omega_1^{\Bb}(\rho)}(\rho,\infty) > \frac{1-\delta(1+\kappa^2L)} {\delta\kappa}\bigg||T_1(\rho)|=\infty \bigg) \\
&\le \frac{\delta\kappa}{1-\delta(1+\kappa^2L)} \E\big[\calR_{\omega_1^{\Bb}(\rho)}(\rho,\infty)\big| |T_1(\rho)|=\infty\big] .
\end{align*}
Conditioned on the survival of the tree $T_1(\rho)$, the backbone tree $T_1^{\Bb}(\rho)$ is a supercritical Galton-Watson tree with conductances bounded by $\kappa$. Hence, \cite[Lemma~9.1]{LyoPemPer96a} yields a finite first moment for the effective resistance above. This implies, using the independence of the indicator functions again,
\begin{align}\label{backbone=1_EW2}
&\E\Bigg[ \sum_{v\in G_L}\mathds 1_{\{P^{v}_{\omega\backslash \omega(v)}(\eta_*=\infty)< \delta\}}\mathds 1_{\{|G_L\cap T_1^{\Bb}|=1\}} \mathds 1_{\{v\notin T_1^{\Bb}\}} \Bigg] \notag\\
&\le C_2 \frac{\delta\kappa}{1-\delta(1+\kappa^2L)} \E\Bigg[ \sum_{v\in G_L} \sum_{u\in G_L\cap T_1(\rho)}
\E\big[ \mathds 1_{\{|T_1(u)|=\infty,\,|T_1(z)|<\infty \text{ for all }z\in G_L\cap T_1(\rho)\backslash\{u\},\,v\notin T_1^{\Bb}\} } \big| G_L \big] \Bigg] \notag \\
&\le C_2 \frac{\delta\kappa}{1-\delta(1+\kappa^2L)}
\E\bigg[\sum_{v\in G_L} \mathds 1_{\{|G_L\cap T_1^{\Bb}|=1,\,v\notin T_1^{\Bb} \} } \bigg]\notag\\
&\le C_2 \frac{\delta\kappa}{1-\delta(1+\kappa^2L)} m^L
\end{align}
for some constant $C_2=C_2(\nu,\alpha,\kappa)>0$. 
Combining \eqref{backbone=1} with the bounds in \eqref{backbone=1_EW1} and \eqref{backbone=1_EW2}, we get
\begin{align}\label{backbone=1_bound}
\E\Bigg[ \sum_{v\in G_L}\mathds 1_{\{P^{v}_{\omega\backslash \omega(v)}(\eta_*=\infty)< \delta\}}\mathds 1_{\{|G_L\cap T_1^{\Bb}|=1\}} \Bigg] 
\le \p\big(|G_L\cap T_1^{\Bb}|=1\big) + C_2 \frac{\delta\kappa}{1-\delta(1+\kappa^2L)} m^L,
\end{align}
which becomes arbitrarily small, when we first chose $L$ large and then $\delta$ small enough.

\subsubsection{The indirect escape probability in the case that the $L$-th generation of $T_1^{\Bb}$ contains more than one vertex}
Lastly, we consider the the expected number of vertices with an indirect escape probability less than $\delta$ on the event, that the $L$-th generation of the backbone tree $T_1^{\Bb}(\rho)$ consists of more than one vertex, which entails that there always exist a path along which the walker can escape no matter of the starting point. In other words, this means that for every vertex $v\in G_L$ there is a vertex $u\neq v$ with $u\in G_L\cap T_1^{\Bb}$. Hence, we get the same lower bound for the indirect escape probability as in \eqref{backbone=1_EW2_indirectescapeprob}, which implies that $\p_{\omega\backslash\omega(v)}^v(\eta_*=\infty)<\delta$ can only hold if
\begin{align*}
\calR_{\omega_1^{\Bb}(u)}(u,\infty) > \frac{1-\delta(1+\kappa^2L)}{\delta\kappa}.
\end{align*}
We denote the vertices in the $L$-th generation of $T$ by $z_1,z_2,\dots$ (ordered as in the Ulam-Harris tree). Then, with the same arguments as before, we obtain
\begin{align}\label{backbone>1_bound}
&\E\Bigg[ \sum_{v\in G_L}\mathds 1_{\{P^{v}_{\omega\backslash \omega(v)}(\eta_*=\infty)< \delta\}}\mathds 1_{\{|G_L\cap T_1^{\Bb}|>1\}} \Bigg] \notag\\
&= \E\Bigg[ \sum_{v\in G_L} 
\mathds 1_{\{P^{v}_{\omega\backslash \omega(v)}(\eta_*=\infty)< \delta\}} 
\mathds 1_{\{|G_L\cap T_1^{\Bb}|>1\}} 
\sum_{i=1}^{|G_L|} \mathds 1_{\{z_i\in T_1(\rho)\cap G_L,\, z_i\neq v,\, |T_1(z_i)|=\infty,\, |T_1(z_j)|<\infty \text{ for all } j<i\}}\Bigg]\notag\\
&\le \E\Bigg[ \sum_{v\in G_L} \sum_{i=1}^{|G_L|} \mathds 1_{\{z_i\in T_1(\rho)\cap G_L,\, z_i\neq v\}}
\E\bigg[ \mathds 1_{\big\{ \calR_{\omega_1^{\Bb}(z_i)}(z_i,\infty) > \frac{1-\delta(1+L)} \delta,\,|T_1(z_i)|=\infty \big\}}  
\mathds 1_{\{|T_1(z_j)|<\infty \text{ for all } j<i\}}\bigg| G_L \bigg]\Bigg]\notag\\
&\le \E\Bigg[ \sum_{v\in G_L} \sum_{i=1}^{|G_L|} \mathds 1_{\{z_i\in T_1(\rho)\cap G_L,\, z_i\neq v\}}
\E\big[\mathds 1_{\{|T_1(z_j)|<\infty \text{ for all } j<i\}}\big| G_L \big]
\E\bigg[\mathds 1_{\big\{ \calR_{\omega_1^{\Bb}(z_i)}(z_i,\infty) > \frac{1-\delta(1+\kappa^2L)}{ \delta\kappa},\,|T_1(z_i)|=\infty \big\}} \bigg] \Bigg] \notag\\
&\le \p\bigg(\calR_{\omega_1^{\Bb}(\rho)}(\rho,\infty) > \frac{1-\delta(1+\kappa^2L)}{\delta\kappa} \bigg| |T_1(\rho)|=\infty \bigg)
\notag\\&\hspace{0.8cm} \times
\E\Bigg[ \sum_{v\in G_L} \sum_{i=1}^{|G_L|} \mathds 1_{\{z_i\in T_1(\rho)\cap G_L,\, z_i\neq v\}}
\E\big[\mathds 1_{\{|T_1(z_j)|<\infty \text{ for all } j<i\}}\big| G_L \big]
\p(|T_1(z_i)|=\infty) \Bigg] \notag\\
&\le C_2 \frac{\delta\kappa}{1-\delta(1+\kappa^2L)} \E\Bigg[ \sum_{v\in G_L} \sum_{i=1}^{|G_L|} \mathds 1_{\{z_i\in T_1(\rho)\cap G_L,\, z_i\neq v,\,|T_1(z_i)|=\infty,\, |T_1(z_j)|<\infty \text{ for all } j<i \}}
\Bigg]\notag\\
&\le C_2 \frac{\delta\kappa}{1-\delta(1+\kappa^2L)} m^L
\end{align}
for some constant $C_2>0$ and $\delta>0$ small enough. For $L$ fixed this bound vanishes as $\delta$ tends to zero.

In total, adding up \eqref{backbone=0_bound}, \eqref{backbone=1_bound} and \eqref{backbone>1_bound} yields the bound
\begin{align*}
\E\Bigg[ \sum_{v\in G_L}\mathds 1_{\{P^{v}_{\omega\backslash \omega(v)}(\eta_*=\infty)< \delta\}} \Bigg]
&\le h(L,\nu,\alpha) + C_1\frac\delta{1-\delta(1+\kappa L)} m^L\\
&+ \p\big(|G_L\cap T_1^{\Bb}|=1\big) + C_2 \frac{\delta\kappa}{1-\delta(1+\kappa^2L)} m^L\\
&+ C_2 \frac{\delta\kappa}{1-\delta(1+\kappa^2L)} m^L\\
&\le 
\gamma(L,\nu,\alpha) + C\frac{\delta\kappa}{1-\delta(1+\kappa^2L)} m^{L},
\end{align*}
where $C=C(\nu,\alpha,\kappa)>0$ is a constant and $\gamma(L,\nu,\alpha)$ converges to zero as $L\to\infty$.
\hfill $\Box$

\subsection{Proof of Lemma~\ref{lem:momentboundregtimes}}
From Lemma~\ref{lem:annealedescapeprob} and Proposition~\ref{prop:stationaryindpendent} follows
\begin{align*}
\bbE\big[(\tau_2-\tau_1)^q\big] 
=\bbE[\tau_1^q|\eta_*=\infty,\xi(\rho,\rho^*)=1]
\le \bbE[\tau_1^q] \bbP(\eta_*=\infty,\xi(\rho,\rho^*)=1)^{-1}
\le C(\varepsilon) \bbE[\tau_1^q] .
\end{align*}
Consequently, the boundedness of the moments of the first regeneration time follows once we show that $P(\tau_1\ge n)\le C(\varepsilon)n^{-r}$ holds for $r$ sufficiently large.

The first regeneration time is large either when the first regeneration occurs in a large generation or when the random walk stays in the first generations for a long time, 
\begin{align}\label{bound-regtimelarge}
\bbP(\tau_1\ge n^3) 
&=\bbP(\tau_1\ge n^3, \tau_1\ge \eta_n) +\bbP(\tau_1\ge n^3, \tau_1< \eta_n)\notag\\
&\le \bbP(\tau_1\ge \eta_n) +\bbP(\eta_n\ge n^3).
\end{align}
Applying Markov's inequality and Lemma~\ref{lem:momentboundregdist} gives an upper bound for the first summand
\begin{align}\label{bound-reginlargegen}
\bbP(\tau_1\ge \eta_n) 
= \bbP(|X_{\tau_1}|\ge n) 
\le \frac{\bbE[|X_{\tau_1}|^r]}{n^r} 
\le C_\varepsilon n^{-r}
\end{align}
for any arbitrary $r>0$ and some constant $C_\varepsilon=C_\varepsilon(r)>0$. 

To derive a polynomially small bound for the second summand of \eqref{bound-regtimelarge} we follow the arguments of \cite{Aidekon2010LargeDeviations}. The hitting time of the $n$-th generation is large either when the random walk visits a great number of distinct vertices or when the walker returns to one vertex many times. To formalize this we let $\pi_k$ be the $k$-th distinct vertex visited by the random walk and recall the definition of the local time $L(z)$ in \eqref{def:localtime}. We now estimate the probability that the walker requires more than $n^3$ steps to reach the $n$-th generation as explained above:
\begin{align}\label{bound_hittingtimelarge}
	\bbP(\eta_n>n^3) 
	\le \bbP(|\{X_0,\dots,X_{\eta_n}\}|> n^2) 
	+ \bbP(\exists k\le n^2: L(\pi_k)> n).
\end{align}
We will establish upper bounds for both probabilities separately.

\subsubsection{Bound on $\bbP(|\{X_0,\dots,X_{\eta_n}\}|> n^2) $}

The random walk can hit more than $n^2$ distinct vertices before hitting the $n$-th generation only if the walker sees more than $n$ distinct vertices in at least one generation. However, this implies that the random walk has to return to the ancestor of each of these vertices. Using the uniform bound in Lemma~\ref{lem:annealedescapeprob} on the probability of hitting the ancestor of the root, we can show that the probability of this decays exponentially in $n$. To specify this, we recall that $G_k$ denotes the $k$-th generation of the tree. Then, we have
\begin{align}\label{bound:distinctverticesSum}
	\bbP(|\{X_0,\dots,X_{\eta_n}\}|> n^2) 
\le \sum_{k=1}^n \bbP(|\{X_0,\dots,X_{\eta_n}\}\cap G_k|\ge n).
\end{align}
We set $t_1^k=\eta_k$ and recursively for $i>1$
\begin{align*}
	t_i^k = \inf\{m>t_{i-1}^k: |X_m|=k,\, X_m\neq X_l \text{ for all }l<m \},
\end{align*}
so that $t_i^k$ is the hitting time of the $i$-th distinct vertex in generation $k$. Then,  
\begin{align*}
	\bbP(|\{X_0,\dots,X_{\eta_n}\}\cap G_k|\ge n)
\le \bbP( t_n^k <\infty) . 
\end{align*}
If $t_n^k$ is finite, $t_{n-1}^k$ has to be finite and the random walk has to return to the parent of vertex $X_{t_{n-1}^k}$, which implies 
\begin{align}\label{bound_distictverticesGenk}
\bbP( t_n^k <\infty) 
	\le  \E\left[ \sum_{v\in T_k} P_\omega^v(\eta_{v^*}<\infty) P_\omega(t_{n-1}^k<\infty, X_{t_{n-1}^k}=v)\mathds 1_{\{v\in T\}}\right] .
\end{align}
The return probability to the ancestor $v^*$ of $v$ will increase when the conductance of the edge $(v^*,v)$ is replaced by $\kappa$. To formalize this, we denote this modification of a given environment $\omega=(T,\rho,\xi)$ by $\tilde\omega_v=(T,\rho,\tilde\xi)$ with
\begin{align*}
	\tilde\xi(z^*,z) = \begin{cases}
		\xi(z^*,z),& z\ne v\\
		\kappa ,& z=v.
	\end{cases}
\end{align*}
We then have 	$P_\omega^v(\eta_{v^*}<\infty)\leq P_{\tilde\omega_v}^v(\eta_{v^*}<\infty)$.
Due to \eqref{bound_distictverticesGenk}, this implies
\begin{align*}
	\bbP( t_n^k <\infty) 
	&\le  \E\left[ \sum_{v\in T_k} P_{\tilde\omega_v}^v(\eta_{v^*}<\infty) P_\omega(t_{n-1}^k<\infty, X_{t_{n-1}^k}=v)|v\in T\right] \\
	&=  \E\left[ \sum_{v\in T_k} P_{\tilde\omega_v}^v(\eta_{v^*}<\infty)|v\in T] \E[P_\omega(t_{n-1}^k<\infty, X_{t_{n-1}^k}=v)\right] \\
	&=\E[P_{\tilde\omega_\rho}^\rho(\eta_*<\infty)] \E\left[ \sum_{v\in T_k} P_\omega(t_{n-1}^k<\infty, X_{t_{n-1}^k}=v)\right] \\
	&=\E[P_{\omega}^\rho(\eta_*<\infty|\xi(\rho^*,\rho)=\kappa)] \bbP(t_{n-1}^k<\infty)\\
	&=\bbP(\eta_*<\infty|\xi(\rho^*,\rho)=\kappa) \bbP(t_{n-1}^k<\infty)
\end{align*}
To see the first equality, observe that the random variables $ P_\omega(t_{n-1}^k<\infty, X_{t_{n-1}^k}=v)$ and $P_{\tilde\omega_v}^v(\eta_{v^*}<\infty)$ are independent under $\p(\,\cdot\,| v\in T)$.
Iterating this gives rise to
\begin{align*}
	\bbP( t_n^k <\infty) \le \bbP(\eta_*<\infty|\xi(\rho^*,\rho)=\kappa)^n. 
\end{align*}
Analogously to Lemma~\ref{lem:annealedescapeprob} we can show
\begin{align*}
\bbP(|\{X_0,\dots,X_{\eta_n}\}\cap G_k|\ge n)
\le	\bbP( t_n^k <\infty)
\le (1-c_\varepsilon)^n
\end{align*}
for some constant $c_\varepsilon=c_\varepsilon(\kappa)>0$. Finally, due to \eqref{bound:distinctverticesSum}, we arrive at 
\begin{align}\label{bound:distinctvertices}
	\bbP(|\{X_0,\dots,X_{\eta_n}\}|> n^2) 
	\le n (1-c_\varepsilon)^n,
\end{align}
which decays exponentially in $n$.

\subsubsection{Bound on $\bbP(\exists k\le n^2: L(\pi_k)> n)$}
Next we will study the probability that the walker returns to one of the vertices $\pi_1,\dots,\pi_{n^2}$ at least $n$ times. Lemma~\ref{lem:returnprob} gives an estimate for the probability of returning to the root at least $n$ times. Using this result we can derive a bound on the probability of visiting the vertex $\pi_k$ at least $n$ times. 

We let $\eta_v^{(1)} = \eta_v = \inf\{m\ge 0 : X_m=v\}$ be the hitting time of the vertex $v\in T$ and $\eta_v^+ = \inf\{m\ge 1: X_m=v\}$ be the return time to vertex $v$. Further, we define recursively for $k>1$
\begin{align*}
	\eta_v^{(k)} = \inf\{ m> \eta_v^{(k-1)} : X_m = v\},
\end{align*}
which indicates when the random walk hits the vertex $v$ the $k$-th time. Now we can write
\begin{align}\label{bound_Prob_LargeNumberOfHittings_kthVertex}
	 \bbP(L(\pi_k)> n) \notag  
	 &=  \E\left[ \sum_{v\in T} P_\omega(\eta_v^{(n+1)}<\infty, \pi_k=v) \right] \notag \\
	&=  \E\left[ \sum_{v\in T} P_\omega(\eta_v<\infty, \pi_k=v) P^v_\omega(\eta_v^+<\infty)^{n} \right] \notag \\
	&\le   \E\left[ \sum_{v\in T} P_\omega(\eta_v<\infty, \pi_k=v) (1- P^v_\omega(\eta_v^+=\infty, \eta_{v^*}=\infty))^{n} \right] .
\end{align}
Recall that $\tilde\omega_v$ is the environment where the conductance of the edge $(v^*,v)$ is replaced by $\kappa$, such that 
\begin{align*}
	P^v_\omega(\eta_v^+=\infty, \eta_{v^*}=\infty) 
	\geq  P^v_{\tilde\omega_v}(\eta_v^+=\infty, \eta_{v^*}=\infty)
\end{align*}
and the lower bound is independent of $P_\omega(\eta_v<\infty, \pi_k=v)$. Together with \eqref{bound_Prob_LargeNumberOfHittings_kthVertex}, we get
\begin{align*}
	\quad&\bbP(L(\pi_k)> n)\\
	&\le  \E\left[ \sum_{v\in\bbT} P_\omega(\eta_v<\infty, \pi_k=v) (1- P^v_{\tilde\omega_v}(\eta_v^+=\infty, \eta_{v^*}=\infty))^{n}\right] \\
	&= \E[(1- P^\rho_{\tilde\omega_\rho}(\eta_\rho^+=\infty, \eta_*=\infty))^{n}] \E\left[ \sum_{v\in T} P_\omega(\eta_v<\infty, \pi_k=v) \right]  \\
	&= \E[(1- P^\rho_{\tilde\omega_\rho}(\eta_\rho^+=\infty, \eta_*=\infty))^{n}]\\
	&= \E[P^\rho_{\tilde\omega_\rho}(\eta_\rho^+<\infty)^{n}].
\end{align*}
By Lemma~\ref{lem:returnprob} this bound decays polynomially in $n$. This implies that for any $\beta>0$ there exists a constant $C_\varepsilon=C_\varepsilon(\beta)$ such that
\begin{align*}
	\bbP(\exists k\le n^2: L(\pi_k)> n)
	\leq \sum_{k=1}^{n^2} \bbP(L(\pi_k)> n) 
	\le C_\varepsilon n^{2-\beta}.
\end{align*}
Combining this with the exponentially small bound in \eqref{bound:distinctvertices} we obtain the polynomial decay of the probability in \eqref{bound_hittingtimelarge}:
\begin{align}
	\bbP(\eta_n>n^3) 
	\le C'_\varepsilon n^{2-\beta}.
\end{align}
Finally, due to \eqref{bound-regtimelarge} and \eqref{bound-reginlargegen}, we arrive at
\begin{align*}
	\bbP(\tau_1\ge n^3) 
	\le Cn^{-r} + C'_\varepsilon n^{2-\beta}.
\end{align*}
Since we can choose $r,\,\beta>0$ arbitrarily, this shows that the $q$-th moment of the first regeneration time $\bbE[\tau_1^q]$ remains bounded, which completes the proof.
\hfill $\Box$

\subsection{Proof of Lemma~\ref{lem:returnprob}}
Recalling that $L(\rho)$ counts the number of visits to the root, we have
\begin{align*}
	\bbP(L(\rho)>n|\xi(\rho^*,\rho)=\kappa) 
= \E[P_{\tilde\omega_\rho}(\eta_\rho^+<\infty)^n].
\end{align*}
Here $\tilde\omega_\rho$ denotes the environment where the conductance of the edge $(\rho^*,\rho)$ is replaced by $\kappa$ and $\eta_\rho^+=\inf\{m\ge 1:X_m=\rho\}$ is the first return time.

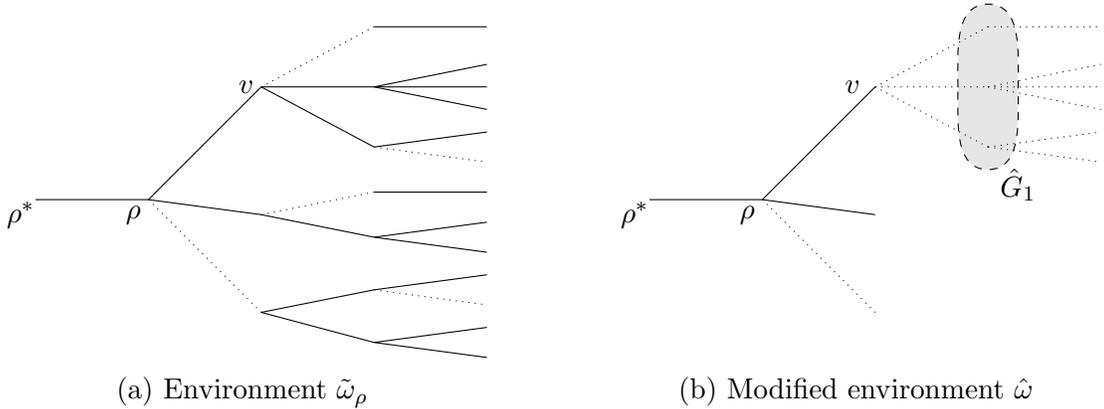
\begin{figure}
	\centering
	\begin{minipage}[t]{.49\linewidth}
	\centering
	\begin{tikzpicture}
		\draw[-] (0,0) to (-1.5,0);
		\draw (-0.2,-0.2) node[] {\small$\rho$};
		\draw (-1.7,-0.2) node[] {\small$\rho^*$};
		\draw[-] (0,0) to (1.5,1.5);
		\draw[-] (0,0) to (1.5,-0.2);
		\draw (1.3,1.5) node[] {\small$v$};
		\draw[style=dotted] (0,0) to (1.5,-1.5);
		\draw[-] (1.5,1.5) to (3,1.5);
		\draw[style=dotted] (1.5,1.5) to (3,2.3);
		\draw[-] (1.5,1.5) to (3,0.7);
		\draw[-] (1.5,-1.5) to (3,-1.2); 
		\draw[-] (1.5,-1.5) to (3,-1.9);
		\draw[style=dotted] (1.5,-0.2) to (3,0.1);
		\draw[-] (1.5,-0.2) to (3,-0.5);
		\draw[-] (3,-1.9) to (4.5,-1.7); 
		\draw[style=dotted] (3,-1.2) to (4.5,-1.4); 
		\draw[-] (3,-1.9) to (4.5,-2.1); 
		\draw[-] (3,-1.2) to (4.5,-1); 
		\draw[-] (3,1.5) to (4.5,1.5); 
		\draw[-] (3,1.5) to (4.5,1.8); 
		\draw[-] (3,1.5) to (4.5,1.2); 
		\draw[-] (3,0.7) to (4.5,0.9); 
		\draw[style=dotted] (3,0.7) to (4.5,0.5); 
		\draw[-] (3,2.3) to (4.5,2.3); 
		\draw[-] (3,0.1) to (4.5,0.1);
		\draw[-] (3,-0.5) to (4.5,-0.3);
		\draw[-] (3,-0.5) to (4.5,-0.7);
	\end{tikzpicture}
	\subcaption{Environment $\tilde\omega_\rho$}
\end{minipage}
\hfill
	\begin{minipage}[t]{.49\linewidth}
	\centering
	\begin{tikzpicture}
		\draw[-] (0,0) to (-1.5,0);
		\draw (-0.2,-0.2) node[] {\small $\rho$};
		\draw (-1.7,-0.2) node[] {\small $\rho^*$};
		\draw[-] (0,0) to (1.5,1.5);
		\draw[-] (0,0) to (1.5,-0.2);
		\draw (1.2,1.5) node[] {\small $v$};
		\draw[style=dotted] (0,0) to (1.5,-1.5);
		\filldraw[fill=gray!20,draw=black,style=dashed] (2.6,1.5) to [out=90, in=180] (3,2.6)to [out=0, in=90] (3.4,1.5) to [out=270, in=0] (3,0.4) to [out=180, in=270] (2.6,1.5);
		\draw (3.4,0.2) node[] {\small $\hat G_1$};
		\draw[style=dotted] (1.5,1.5) to (3,1.5);
		\draw[style=dotted] (1.5,1.5) to (3,2.3);
		\draw[style=dotted] (1.5,1.5) to (3,0.7);
		\phantom{\draw[-] (3,-1.9) to (4.5,-2.1);} 
		\draw[style=dotted] (3,1.5) to (4.5,1.5); 
		\draw[style=dotted] (3,1.5) to (4.5,1.8); 
		\draw[style=dotted] (3,1.5) to (4.5,1.2); 
		\draw[style=dotted] (3,0.7) to (4.5,0.9); 
		\draw[style=dotted] (3,0.7) to (4.5,0.5); 
		\draw[style=dotted] (3,2.3) to (4.5,2.3); 
	\end{tikzpicture}
	\subcaption{Modified environment $\hat\omega$}
\end{minipage}
\caption{$v$ is the first descendant of $\rho$. In the modified environment $\hat\omega$ the conductance of each edge in $T(v)$ is replaced by $\varepsilon$ and all edges in $T(z)$ for $z\in G_1(T)$, $z\ne v$ are removed. The first generation $\hat G_1$ of the subtree $T(v)$ is marked gray. Edges with conductance $\varepsilon$ are indicated by dotted lines;  edges with larger conductance are indicated by solid lines.}
	\label{fig:modifiedtree_hat}
\end{figure}

We will start by determining a bound on the quenched return probability $P_{\tilde\omega_\rho}(\eta_\rho^+<\infty)$ by modifying the environment once more. We let $v$ be the first descendant of the root and we replace the conductances of all edges in the subtree $T(v)$ by $\varepsilon$. The first generation of this subtree is denoted by $\hat G_1=G_1(T(v))$. Further, we remove all edges in the other subtrees $T(z)$ rooted at a descendant $z$ of $\rho$ with $z\neq v$. We call the new environment $\hat\omega=(T,\rho, \hat\xi)$ with 
\begin{align*}
	\hat\xi(z^*,z) =\begin{cases}
		\kappa,& z=\rho\\
		\xi(z^*,z), & |z|=1\\
		\varepsilon,& z,\,z^*\in T(v)\\
		0,& \text{otherwise,}
	\end{cases}
\end{align*}
see Figure~\ref{fig:modifiedtree_hat} for an example. 

A straightforward conductance argument shows that 
\begin{align*}
	P_{\tilde\omega_\rho}^\rho(\eta_{\rho}^+<\infty) 
\leq  P_{\hat\omega}^\rho(\eta_\rho^+<\infty), 
\end{align*}
which implies 
\begin{align}\label{bound_nreturns_newenvironment}
	\E[P_{\tilde\omega_\rho}^\rho(\eta_{\rho}^+<\infty)^{n}]
	\le \E[P_{\hat\omega}^\rho(\eta_{\rho}^+<\infty)^{n}]
	= \E[P_{\hat\omega}(L(\rho)> n)].
\end{align}
We introduce
\begin{align*}
	H_n = \sum_{k=2}^{\eta_\rho^{(n+1)}} \mathds 1_{\{X_{k-2}=\rho,X_{k-1}=v,X_k\in \hat G_1\}},
\end{align*}
that indicates how often the random walk moves from $\rho$ to a vertex in $\hat G_1$ before it returns to $\rho$ for the $n$-th time.
Then, we can write 
\begin{align}\label{nReturnsToRho}
	\E[P_{\hat\omega}(L(\rho)> n)]
	= \E[P_{\hat\omega}(L(\rho)> n,H_n\ge n\delta_n)] +\E[P_{\hat\omega}(L(\rho)> n,H_n< n\delta_n)],
\end{align}
where we choose $\delta_n= c_\varepsilon n^{-\frac13}$ with $c_\varepsilon = \frac{\varepsilon^2}{\kappa^2+\varepsilon\kappa}$.
We will treat the two summands separately.

\subsubsection{Bound on $\E[P_{\hat\omega}(L(\rho)> n,H_n> n\delta_n)]$}

For this bound we compare the return probabilities of the random walk on $\hat\omega$ with those of a simple random walk on a Galton-Watson tree with unit conductance.

In the environment $\hat\omega$ each edge in the subtree $T(v)$ has the same conductance. For this reason the random walk, observed only at the vertices of $T(v)$, has the same distribution as a simple random walk on $T(v)$. 
On the event $\{L(\rho)> n,H_n> n\delta_n\}$ the random walk on $T$ has to return from $\hat G_1$ to $\rho$ at least $n\delta_n$ times. In particular, this means that the walk restricted to $T(v)$ has to return from $\hat G_1$ to $v$ at least $n\delta_n$ times. Therefore, it suffices to bound the probability that a simple random walk on a Galton-Watson tree hits the root at least $n\delta_n$ times. 
As proven by \cite{Piau98CLT}, the probability that the first regeneration time of a random walk on a Galton-Watson tree with unit conductances and without leaves is at least $n$ is bounded by $\exp(-c'n^{1/6})$ for some constant $c'>0$. This implies the following stretched exponential bound:
\begin{align}\label{bound_LG1>ndelta}
	\E[P_{\hat\omega}(L(\rho)> n,H_n> n\delta_n)] 
	\le \mathrm e^{-c'  (\lfloor n\delta_n\rfloor )^{\frac16} }
	\le \mathrm e^{-\tilde c_\varepsilon n^{\frac19} }
\end{align}
for some constant $\tilde c_\varepsilon>0$.

\subsubsection{Bound on $\E[P_{\hat\omega}(L(\rho)> n,H_n\le n\delta_n)]$}
We introduce
\begin{align*}
	R_n=n-H_n=\sum_{k=2}^{\eta_\rho^{(n+1)}} \mathds 1_{\{X_{k-2}=X_k=\rho\}} 
	= \sum_{k=1}^n \mathds 1_{\{ \eta_\rho^{(k+1)} - \eta_\rho^{(k)} =2\}}, 
\end{align*}
which is well-defined on the event $\{L(\rho)> n\}=\{\eta_\rho^{(n+1)}<\infty\}$. That is, $R_n$ counts the number of the excursions of length $2$, starting from the root, until the random walk hits the root for the $(n+1)$-th time. Hence, we can rewrite the second summand in \eqref{nReturnsToRho}
\begin{align*}
	E[P_{\hat\omega}(L(\rho)> n,H_n\le n\delta_n)] = 	E[P_{\hat\omega}(L(\rho)> n,R_n\ge n(1-\delta_n))].
\end{align*}
Let $\omega=(T,\rho,\xi)$ be an environment and $\hat\omega$ is the modification introduced above. Under $P_{\hat\omega}^\rho(\,\cdot\,|L(\rho)>n)$, the random variables $\eta_{\rho}^{(2)}-\eta_{\rho}^{(1)},\dots,\eta_{\rho}^{(n+1)}-\eta_{\rho}^{(n)}$ are independent and identically distributed. Hence, as a sum of i.i.d. Bernoulli random variables $R_n$ is binomially distributed with parameters $n$ and success probability 
\begin{align*}
	P_{\hat\omega}^\rho(\eta_{\rho}^{(k)}-\eta_{\rho}^{(k-1)}=2|L(\rho)>n)= P_{\hat\omega}^\rho(\eta_\rho^+=2|\eta_\rho^+<\infty)
	= \frac{P_{\hat\omega}^\rho(\eta_\rho^+=2)} {P_{\hat\omega}^\rho(\eta_\rho^+<\infty)}.
\end{align*}
It follows
\begin{align*}
	P_{\hat\omega}(L(\rho)> n,\,R_n\ge n(1-\delta_n)) \le P(B_n\ge n(1-\delta_n)),
\end{align*}
where $B_n$ is a binomially distributed random variable with parameters $n$ and $P_{\hat\omega}^\rho(\eta_\rho^+=2)$.
This probability is given by 
\begin{align*}
	P_{\hat\omega}(\eta_\rho^+=2)
	&=1-P_{\hat\omega}(X_2\ne \rho)
	=1-P_{\hat\omega}(X_1=v)P_{\hat\omega}^v(X_1\in \hat G_1)\\
	&=1-\frac{\xi(\rho,v)}{\kappa+\sum_{z\sim,z\ne\rho^*}\xi(\rho,z)}\cdot\frac{\varepsilon|\hat G_1|}{\xi(\rho,v)+\varepsilon|\hat G_1|}\\
	&\le 1-\frac{\varepsilon}{\kappa|G_1|+\kappa}\cdot\frac{\varepsilon}{\kappa+\varepsilon} \le 1-\hat p
\end{align*}
with
\begin{align*}
	\hat p = \frac{c_\varepsilon}{2|G_1|}, \quad
	c_\varepsilon = \frac{\varepsilon^2}{\kappa^2+\varepsilon\kappa}.
\end{align*}
This implies 
\begin{align}\label{binomialbound}
	P_{\hat\omega}(L(\rho)> n,R_n\ge n(1-\delta_n))
	\le P(\tilde B_n\ge n(1-\delta_n))
	= P(\tilde B_n-n(1-\hat p)\ge n(\hat p-\delta_n)),
\end{align}
where $\tilde B_n$ is a binomially distributed random variable with parameters $n$ and $1-\hat p$. We can use Hoeffding's inequality to establish a stretched exponential bound for the probability above, if $\hat p-\delta_n>0$ holds. To guarantee that $\hat p-\delta_n$ is positive, we introduce the event
\begin{align*}
	D_n = \Big\{|G_1|\le \frac13 n^{\frac13}\Big\}.
\end{align*}
Then, for $\omega\in D_n$ (or equivalently $\hat\omega\in D_n$) we have
\begin{align*}
	\hat p-\delta_n = \frac{c_\varepsilon}{2|G_1|}-c_\varepsilon n^{-\frac13}
	\ge \frac{c_\varepsilon}2 n^{-\frac13}>0.
\end{align*}
Plugging this in \eqref{binomialbound} and applying Hoeffding's inequality afterwards we arrive at
\begin{align*}
	P_{\hat\omega}(L(\rho)> n,R_n\ge n(1-\delta_n))
	\le P\Big(B_n-n(1-\hat p)\ge \frac{c_\varepsilon}2 n^{\frac23}\Big)
	\le \mathrm e^{-\frac{c_\varepsilon^2}2n^{\frac13}}.
\end{align*}
Since this bound holds for all $\omega\in D_n$, we finally get
\begin{align*}
	\E[P_{\hat\omega}(L(\rho)> n,H_n\le n\delta_n)\mathds 1_{D_n}] 
	\le \mathrm e^{-\frac{c_\varepsilon^2}2n^{\frac13}}.
\end{align*}
Further, using the Markov inequality we obtain
\begin{align*}
	\E[P_{\hat\omega}(L(\rho)> n,H_n\le n\delta_n)\mathds 1_{D_n^c}]
	\le P \Big(|G_1|> \frac13 n^{\frac13}\Big)
	\le 3^m E[|G_1|^m]n^{-\frac m3}
\end{align*}
for $m>0$ arbitrary. 
From the last two estimates follows that
\begin{align*}
	\E[P_{\hat\omega}(L(\rho)> n,H_n\le n\delta_n)] \le C_\varepsilon n^{-\frac m3}
\end{align*}
holds for some constant $C_\varepsilon=C_\varepsilon(m)>0$. Together with \eqref{bound_nreturns_newenvironment} and the stretched exponential bound in \eqref{bound_LG1>ndelta} this implies that for any $\beta>0$ there exists some constant $C'_\varepsilon=C'_\varepsilon(\beta)>0$ such that
\begin{align*}
	\E[P_{\tilde\omega_\rho}^\rho(\eta_{\rho}^+<\infty)^{n}]\le
	\E[P_{\hat\omega}(L(\rho)> n)]
	\le C'_\varepsilon n^{-\beta}.
\end{align*}
\hfill $\Box$

\section{Proof of Theorem~\ref{thm:CLT}}

The arguments to obtain a functional CLT from a sufficient regeneration structure have been applied in several varying situations. We may follow the arguments of \cite{Sznitman2000CLT} and only give the main steps in the proof. 

In a first step, consider the jump process with increments
\begin{align*}
Z_k = |X_{\tau_{k}}|- |X_{\tau_{k-1}}|-(\tau_{k}-\tau_{k-1})v, \quad\text{for } k\ge 1
\end{align*}
where
\begin{align*}
v = \frac{\bbE\big[|X_{\tau_2}|- |X_{\tau_1}|\big]}{\bbE[\tau_2-\tau_1]} 
\end{align*}
is the speed of the random walk.

By Proposition~\ref{prop:stationaryindpendent} the random variables $Z_2,Z_3,\dots$ are independent and identically distributed and by Lemma~\ref{lem:momentboundregdist} and Lemma~\ref{lem:momentboundregtimes} they are square integrable with $\bbE[Z_2] = 0$ and $\var(Z_2) = \bbE[Z_2^2] \eqqcolon \sigma^2$. 
Hence, we may apply Donsker's Theorem (see for example~\cite[Theorem~14.1]{billingsley2013convergence}), which implies  
\begin{align}\label{convergence_donsker}
\bigg( \frac1{\sqrt n} \sum_{k=1}^{\lfloor nt\rfloor} Z_k \bigg)_{t\in [0,1]} 
\xrightarrow[n\to\infty]{d} (\sigma B_t)_{t\in [0,1]}.
\end{align}

Next, we correct the jump times of the process considered above by applying a random time change. 
For this purpose we let $k_n$ be the number of regenerations until time $n$, i.e. $k_n$ is the integer, such that $\tau_{k_n}\le n<\tau_{k_n+1}$. Since $\tau_k-\tau_{k-1}, k\geq 2$ are i.i.d. and integrable, 
\begin{align}\label{limit_kn/n}
\frac {k_n}n \xrightarrow{n\to\infty} \frac1{\bbE[\tau_2-\tau_1]} \quad \mathbb{P}-\text{a.s.}
\end{align}
Using a counterpart of Dini's theorem (see \cite[Part~II, Problem~127]{polya1972analysisI}), the following uniform convergence holds: 
\begin{align*}
\sup_{t\in[0,1]}\bigg|\frac{k_{\lfloor nt\rfloor}}n-\frac{t}{\bbE[\tau_2-\tau_1]} \bigg|  \xrightarrow{n\to\infty} 0 \quad \mathbb{P}-\text{a.s.}
\end{align*} 
We may use the arguments of \cite[Section~14]{billingsley2013convergence} to conclude
\begin{align*}
\Big( \frac1{\sqrt n} \sum_{k=1}^{k_{\lfloor nt\rfloor}} Z_k \Big)_{t\in [0,1]}\xrightarrow[n\to\infty]{d} 
\Big(\sigma B_{\frac t{\bbE[\tau_2-\tau1]}}\Big)_{t\in[0,1]}.
\end{align*}
Due to the scale invariance of the Brownian motion we arrive at
\begin{align}\label{convergence_timescaledprocess}
\Big( \frac1{\sqrt n} \sum_{k=1}^{k_{\lfloor nt\rfloor}} Z_k \Big)_{t\in [0,1]}\xrightarrow[n\to\infty]{d} (\tilde\sigma B_t)_{t\in[0,1]}
\end{align}
with $\tilde \sigma = \sigma\bbE[\tau_2-\tau_1]^{-\frac12}$.

Is remains to add the path between regenerations. Let $d_S$ denote the Skorokhod metric, then 
\begin{align}\label{bound-Skorokhoddist}
&d_S\bigg( \Big( \frac1{\sqrt n} \big(|X_{\lfloor nt\rfloor}|-\lfloor nt\rfloor v\big)\Big)_{t\in[0,1]}, 
\Big( \frac1{\sqrt n} \sum_{k=1}^{k_{\lfloor nt\rfloor}} Z_k \Big)_{t\in [0,1]} \bigg)\notag \\
&\qquad\le \sup_{t\in[0,1]} \bigg|  \frac1{\sqrt n} \big(|X_{\lfloor nt\rfloor}|-\lfloor nt\rfloor v\big)- \frac1{\sqrt n} \sum_{k=1}^{k_{\lfloor nt\rfloor}} Z_k \bigg|\notag \\
&\qquad\le \frac1{\sqrt n} \sup_{t\in[0,1]} \big(
\big| |X_{\lfloor nt\rfloor}|-|X_{\tau_{k_{\lfloor nt \rfloor}}}|\big| + \big| \lfloor nt\rfloor - \tau_{k_{\lfloor nt\rfloor}}\big|v \big) \notag 
\\
&\qquad\le \frac{1+v}{\sqrt n}  \sup_{t\in[0,1]}
\big( \lfloor nt\rfloor - \tau_{k_{\lfloor nt\rfloor}} \big)\notag \\
& \qquad \le(1+v) \bigg(\frac{|\tau_1|}{\sqrt n} + \frac{\max_{k\in\{1,\dots,n\}}|\tau_{k+1}-\tau_k|}{\sqrt n}\bigg).
\end{align}
The first summand vanishes in probability, since $\tau_1$ is $\bbP$-almost surely finite. Furthermore, the random variables $|\tau_{k+1}-\tau_k|$, $k=1,\dots,n$, are independent and identically distributed, non-negative and square integrable by Proposition~\ref{prop:stationaryindpendent} and Lemma~\ref{lem:momentboundregtimes}. This implies $\max_{k\in\{1,\dots,n\}}|\tau_{k+1}-\tau_k| = o_{\bbP}(\sqrt{n})$ and then  
\eqref{bound-Skorokhoddist} vanishes in probability.
Finally, due to the convergence result in \eqref{convergence_timescaledprocess} and \cite[Theorem~3.1]{billingsley2013convergence} we get
\begin{align*}
\Big(\frac 1{\sqrt n} \big(|X_{\lfloor nt\rfloor}|-\lfloor nt\rfloor v\big)\Big)_{t\in[0,1]} 
\xrightarrow[n\to\infty]{d} (\tilde\sigma B_t)_{t\in[0,1]}.
\end{align*}
\hfill $\Box$

\section{The volatility for $\varepsilon$ close to 0}

In this section we prove Theorem \ref{thm:volatility}. To emphasize the dependency on $\varepsilon$, we write $\mathrm{P}_\varepsilon$ and $\bbP_\varepsilon$ for the environment law and annealed law, respectively, when the law of the conductances is $\mu$ as in \eqref{def:conductancelaw}. The corresponding expectations are denoted by $\E_\varepsilon$ and $\bbE_\varepsilon$. 
From the proof of Theorem~\ref{thm:CLT} we know that the volatility is given by
\begin{align*}
	\sigma^2(\varepsilon) 
	= \frac{\bbE_\varepsilon[(|X_{\tau_2}|-|X_{\tau_1}| - (\tau_2-\tau_1)v(\mu_\varepsilon))^2]} {\bbE_\varepsilon[\tau_2-\tau_1]},
\end{align*}
where 
\begin{align*}
	v(\mu_\varepsilon)
	= \frac{\bbE_\varepsilon[|X_{\tau_2}|-|X_{\tau_1}|]}
	{\bbE_\varepsilon[\tau_2-\tau_1]}
\end{align*}
is the speed of the random walk. The proof of Theorem~\ref{thm:volatility} follows from the next two lemmas.

\begin{lem}\label{lem:numeratorVolatility}
	There exists a constant $c>0$, independent of $\varepsilon$, such that
	\begin{align*}
		\bbE_\varepsilon[(|X_{\tau_2}|-|X_{\tau_1}|-(\tau_2-\tau_1)v(\mu_\varepsilon))^2]
		\ge c>0.
	\end{align*}
\end{lem}

\begin{lem}\label{lem:denominatorVolatility}
	If $P(|T_1|=\infty)>0$, we have
	\begin{align*}
		\limsup_{\varepsilon\to 0} \bbE_\varepsilon[\tau_2-\tau_1] <\infty.
	\end{align*}
\end{lem}

Combining the representation of the volatility from above with Lemma~\ref{lem:numeratorVolatility} and Lemma~\ref{lem:denominatorVolatility}, we obtain
\begin{align*}
	\liminf_{\varepsilon\to 0} \sigma^2(\varepsilon) 
	\ge \liminf_{\varepsilon\to 0} \frac c{\bbE_\varepsilon[\tau_2-\tau_1]}
	= \frac c{\limsup_{\varepsilon\to 0} \bbE_\varepsilon[\tau_2-\tau_1]}>0.
\end{align*}
\hfill $\Box$

\subsection{Proof of Lemma~\ref{lem:numeratorVolatility}}
We will start by introducing some notations. 
We define $d=\inf\{k\ge 1: p_k>0\}$, which ensures that the probability of a vertex having $d$ descendants is positive. 
Thinking of a tree as a subset of the Ulam-Harris-tree, we let $v_1$ be the first descendant of the root $\rho$ and $v_2$ is the first descendant of $v_1$. 
Now, we may introduce the set of environments $A=A_1\cap A_2$ with
\begin{align*}
	A_1 &= \{\omega\in\Omega: \deg(v) = d+1 \text{ for all } v\in  \{\rho\}\cup G_{1}(T)\},\\
	A_2 &= \{\omega\in\Omega: \xi(v^*,v) = 1 \text{ for all } v\in  \{\rho\}\cup G_{1}(T)\cup G_2(T)\},
\end{align*}
and the set of trajectories
\begin{align*}
	B = \{ X_0=\rho,\, X_1=v_1,\, X_2=v_2\}\cap\{ |X_n|\ge 2 \,\forall n\ge 2\}.
\end{align*}
Then, for all $(\omega,(X_n)_{n\ge 0})\in A\times B_n$ we have
\begin{align*}
	&\tau_1=1,\quad 
	|X_{\tau_{1}}|=|v_1|=1,\quad 
	\tau_2 = 2,\quad 
	|X_{\tau_{2}}|=|v_2|=2.
\end{align*}
This implies
\begin{align}\label{lowerbound_Z2}
	\bbE_\varepsilon[(|X_{\tau_2}|-|X_{\tau_1}|-(\tau_2-\tau_1)v(\mu_\varepsilon))^2]
	&\ge \bbE_\varepsilon[(|X_{\tau_2}|-|X_{\tau_1}|-(\tau_2-\tau_1)v(\mu_\varepsilon))^2 \mathds 1_{A\times B}] \notag \\
	&= (1-v(\mu_\varepsilon))^2 \bbP_\varepsilon(A\times B) \notag \\
	&= (1-v(\mu_\varepsilon))^2 \E_\varepsilon [P_\omega(B)\mathds 1_A(\omega)].
\end{align}
It follows by a straighforward computation and by Lemma~\ref{lem:annealedescapeprob} that there exists a constant $c=c(\nu,\alpha)>0$, independent of $\varepsilon$, such that
\begin{align*}
	\E_\varepsilon [P_\omega(B_n)\mathds 1_A(\omega)] \ge c.
\end{align*}
In order to get a uniform bound for $v(\mu_\varepsilon)$, we compare the speed of the random walk on Galton-Watson trees with random conductances with the speed of a simple random walk on the same tree. \cite{gantert2012random} showed in Theorem~4.4 that random conductances can only slow down the random walk.
Together with the explicit formula for the speed of the simple random walk in Theorem~3.2 of \cite{LyoPemPer96a} we obtain
\begin{align*}
	v(\mu_\varepsilon) \leq  \sum_{k=1}^\infty p_k \frac{k-1}{k+1}<1, 
\end{align*}
We have then a uniform lower bound for \eqref{lowerbound_Z2}, which completes the proof.
\hfill $\Box$

\subsection{Proof of Lemma~\ref{lem:denominatorVolatility}}
In \cite{glatzel2021speed} we studied how the speed $v=v(\mu_\varepsilon)$ of the random walk depends on $\varepsilon$. If $P(|T_1|=\infty)>0$,  Theorem~1.2 of \cite{glatzel2021speed} shows that the limit of the speed for conductances approaching zero is 
\begin{align*}
	\lim_{\varepsilon\to 0} v(\mu_\varepsilon)  = \beta v(\mu_0),
\end{align*}
for some constant $\beta\in (0,1)$. Here $v(\mu_0)$ is the speed of the random walk as usually defined on trees with positive extinction probability, where we condition on $|T_1(\rho)|=\infty$. A formula for $v(\mu_0)$ is given in Remark~4.1 in \cite{gantert2012random}. Using the same arguments as in Theorem~4.2 therein, which states that the speed on weighted Galton-Watson tree without leaves is strictly positive, one can show $v(\mu_0)>0$. Consequently, we have
\begin{align*}
	\lim_{\varepsilon\to 0} v(\mu_\varepsilon)
	= \lim_{\varepsilon\to 0} \frac{\bbE_\varepsilon[|X_{\tau_2}|-|X_{\tau_1}|]} {\bbE_\varepsilon[\tau_2-\tau_1]}>0.
\end{align*}
Due to Lemma~\ref{lem:momentboundregdist}, we can uniformly bound the numerator, which implies
\begin{align*}
	0<\lim_{\varepsilon\to 0} v(\mu_\varepsilon)
		\le  \liminf_{\varepsilon\to 0} \frac C{\bbE_\varepsilon[\tau_2-\tau_1]}
		= \frac C{\limsup_{\varepsilon\to 0} \bbE_\varepsilon[\tau_2-\tau_1]}
\end{align*}
for some constant $C>0$. This can only hold, if
\begin{align*}
	\limsup_{\varepsilon\to 0} \bbE_\varepsilon[\tau_2-\tau_1]<\infty,
\end{align*}
which completes the proof.
\hfill $\Box$

\newpage
\bibliographystyle{alpha}
\bibliography{bibGWT}

\bigskip

{\footnotesize
\noindent
TU Dortmund, \\
Fakult\"at f\"ur Mathematik, \\
Vogelpothsweg 87, 
44227 Dortmund, 
Germany, \\
tabea.glatzel@tu-dortmund.de\\
jan.nagel@tu-dortmund.de
}

\end{document}